\documentclass[10pt]{amsart}
\usepackage[T1]{fontenc}
\usepackage{geometry}
\usepackage[latin1] {inputenc}
\usepackage{amsmath}
\usepackage{amsfonts, amssymb, textcomp}
\usepackage[colorlinks=flase, linkcolor=red,urlcolor=green, citecolor=blue]{hyperref}
\usepackage{subeqnarray}

\usepackage{latexsym}
\usepackage{fancyhdr}
\usepackage{longtable}
\usepackage{amsmath, amssymb}
\usepackage{graphicx}
\usepackage{enumerate}

\numberwithin{equation}{section}

\theoremstyle{plain}
\newtheorem{proposition}{Proposition}[section]
\newtheorem{theorem}[proposition]{Theorem}
\newtheorem{lemma}[proposition]{Lemma}
\newtheorem{corollary}[proposition]{Corollary}
\newtheorem{definition}[proposition]{Definition}

\newtheorem{remark}[proposition]{Remark}

\newcommand{\RR}{\mathbb{R}}
\newcommand{\CC}{\mathbb{C}}
\newcommand{\NN}{\mathbb{N}}

\newcommand{\id}{\operatorname{id}}

\let\on=\operatorname

\newenvironment{demo}[1]{\par\smallskip\noindent{\bf #1.}}{\par\smallskip}

\makeatletter
\@namedef{subjclassname@2020}{%
  \textup{2020} Mathematics Subject Classification}
\makeatother

\newsavebox{\fmbox}
\newenvironment{fmpage}[1]
 {\begin{lrbox}{\fmbox}\begin{minipage}{#1}}
 {\end{minipage}\end{lrbox}\fbox{\usebox{\fmbox}}}

\title[On inclusion relations between weighted spaces of entire functions]
{On inclusion relations between weighted spaces of entire functions}

\author[G.~Schindl]{Gerhard Schindl}

\address{G.~Schindl: Fakult\"at f\"ur Mathematik, Universit\"at Wien, Oskar-Morgenstern-Platz~1, A-1090 Wien, Austria.}
\email{gerhard.schindl@univie.ac.at}

\begin{document}

\begin{abstract}
We characterize the inclusions of weighted classes of entire functions in terms of the defining weights resp. weight systems. First we treat weights defined in terms of a so-called associated weight function where the weight(system) is based on a given sequence. The abstract weight function case is then reduced to the weight sequence setting by using the so-called associated weight sequence. As an application of the main statements we characterize closedness under point-wise multiplication of these classes.
\end{abstract}

\thanks{G. Schindl is supported by FWF-Project 10.55776/P33417}
\date{\today}

\maketitle


\section{Introduction}\label{Introduction}
Weighted spaces of entire functions are defined as follows:
\begin{equation}\label{weightedclassdef}
H^{\infty}_v(\CC):=\{f\in H(\CC): \|f\|_v:=\sup_{z\in\CC}|f(z)|v(|z|)<+\infty\}.
\end{equation}
Here $H(\CC)$ denotes the class of entire functions and $v:[0,+\infty)\rightarrow(0,+\infty)$ is the (radial and strictly positive) {\itshape weight function.} From now on, except stated explicitly otherwise, we assume that $v$ is
\begin{itemize}
	\item[$(*)$] continuous,
	
	\item[$(*)$] non-increasing and
	
	\item[$(*)$] rapidly decreasing, i.e. $\lim_{t\rightarrow+\infty}t^kv(t)=0$ for all $k\ge 0$.
\end{itemize}
This is the same notion to be a weight (function) as it has been considered in \cite{BonetTaskinen18}, \cite{BonetLuskyTaskinen19}, \cite{solidassociatedweight} and in \cite{Bonet2022survey}. In \cite[Sect. 2]{Bonet2022survey} weighted classes are even introduced for weights which are not necessarily non-increasing and rapidly decreasing; such $v$'s are called {\itshape general weights}. Also in \cite{BierstedtBonetTaskinen98}, \cite{BonetDomanskiLindstroemTaskinen98} more general notions have been considered. It is known that $H^{\infty}_v(\CC)$ is a Banach space w.r.t. the norm $\|\cdot\|_v$ and similarly such weighted classes can be defined when $\CC$ is replaced by any other (non-empty) open subset $G\subseteq\CC$ and $v: G\rightarrow(0,+\infty)$. Another very prominent example in the literature is the unit disc $G=\mathbb{D}$. In this article we will focus on $G=\CC$ except in Section \ref{essentialweightsect} and so, in order to lighten notation, we frequently write $H^{\infty}_v$ instead of $H^{\infty}_v(\CC)$ when the case is clear (except e.g. in Section \ref{essentialweightsect}). In Section \ref{weightsection} all relevant information on weights, sequences, weight systems and the corresponding weighted classes is gathered.

Such weighted spaces (and operators on them) have been studied during the last years from different aspects. Very recently in \cite{Bonet2022survey} a detailed survey of available results on this topic is given; see also the numerous literature citations therein.\vspace{6pt}

The aim of this paper is to study and to compare growth relations between weights and weight systems and inclusion relations of the corresponding weighted spaces. To be more precise, let $v$ and $w$ be two weights and write $v\hypertarget{ompreceq}{\preceq}w$ if
\begin{equation}\label{bigOrelation}
w(t)=O(v(t))\;\text{as}\;t\rightarrow+\infty.	
\end{equation}

Let us call $v$ and $w$ {\itshape equivalent,} written $v\hypertarget{sim}{\sim}w$, if
$$v\hyperlink{ompreceq}{\preceq}w\;\text{and}\;w\hyperlink{ompreceq}{\preceq}v.$$
Note that relation \hyperlink{sim}{$\sim$} is denoted by the symbol $\approx$ in \cite{Bonet2022survey} and by $\asymp$ in \cite{AbakumovDoubtsov15} and \cite{AbakumovDoubtsov18}. Obviously, by definition one has
\begin{equation}\label{principalequ}
v\hyperlink{ompreceq}{\preceq}w\Longrightarrow H^{\infty}_{v}(\CC)\subseteq H^{\infty}_{w}(\CC),
\end{equation}
with continuous inclusion (and similarly if $\CC$ is replaced by any other subset $G\subseteq\CC$). So equivalent weights yield equivalent weighted classes (as locally convex vector spaces) and the main question is: What can be said about the implication $\Longleftarrow$ in \eqref{principalequ}, i.e., can the inclusion (resp. equality) of weighted spaces be characterized in terms of a growth relation between two given weights?\vspace{6pt}

In Section \ref{essentialweightsect} we are concerned with this problem. As the anonymous referee has commented in his report, in fact the desired characterization follows from a more general result in \cite{BonetDomanskiLindstroemTaskinen98} dealing with the composition operator $C_{\varphi}(f):=f\circ\varphi$ for $G=\mathbb{D}$ and it also holds for $G=\CC$; see Theorem \ref{essentialweightcharactnew}. However, the obtained characterization does not involve the weights $v$ resp. $w$ directly but a technical modification; see \eqref{assoweightw} and \eqref{assoweightv}. In order to get now the desired equivalence \eqref{principalequ} (for $G=\mathbb{D}$ and $G=\CC$) one has to assume that $v$ is {\itshape essential.}

We also give a second independent proof of this statement: If $v$ is {\itshape essential,} then alternatively by using a consequence of the recent works \cite{AbakumovDoubtsov15} and \cite{AbakumovDoubtsov18} in \eqref{principalequ} equivalence holds true (again for $G=\mathbb{D}$ and $G=\CC$); see Theorem \ref{essentialweightcharact}.\vspace{6pt}

In this work we are basically interested in the case when $M\in\RR_{>0}^{\NN}$ is a given (weight) sequence which satisfies standard growth assumptions and
$$v=v_M: t\mapsto\exp(-\omega_M(t)).$$
Here, $\omega_M$ denotes the so-called {\itshape associated weight function}; see Section \ref{assofctsect}. Such particular weights have been considered in \cite{solidassociatedweight} where solid hulls and cores of the corresponding weighted classes have been studied. The advantage in this setting is that $v_M$ is based on a given sequence and so desired properties for $v_M$ can be expressed and characterized in terms of growth and regularity assumptions for $M$ (known and used in the ultradifferentiable and ultraholomorphic setting). Therefore, more structure and information is involved and available. For example, by combining another recent result from \cite{AbakumovDoubtsov18} and the crucial Lemma \ref{monomialnormlemma}, we prove in Theorem \ref{essentialweightsequthm} that any $v_M$ is essential (for the weighted entire case) and hence in Theorem \ref{weightholombysequcharactsingle} the equivalence $v_M\hyperlink{ompreceq}{\preceq}v_N\Longleftrightarrow H^{\infty}_{v_M}(\CC)\subseteq H^{\infty}_{v_N}(\CC)$ is characterized by the following growth relation between $M$ and $N$:
\begin{equation}\label{strongNMrelation}
\exists\;A\ge 1\;\forall\;j\in\NN:\;\;\;N_j\le AM_j.
\end{equation}

When an abstract weight $v$ is given, the idea is to reduce this case to the weight sequence setting. This can be achieved by considering the so-called {\itshape associated weight sequence} $M^v$; see Section \ref{weighttosequsect}. Note that in view of Remarks \ref{convnoneed}, \ref{assoweightsequ} in order to involve $M^v$ our standard assumptions on $v$ to be a weight are convenient. Moreover, it turns out that in order to use $M^v$ it suffices to assume that $v$ is \emph{(log)-convex} which is in general weaker than being essential in the weighted entire case; see Remark \ref{AtoEremark}. The comments summarized in this Remark together with Theorem \ref{assoweightfctprop} give also a comparison between essential and (log)-convex weights which satisfy the crucial technical assumption {\itshape moderate growth;} see Definition \ref{admissdef1}.\vspace{6pt}

\eqref{strongNMrelation} is much more restrictive than the (known) comparison relation
\begin{equation}\label{weightsequencrelation}
\sup_{j\in\NN_{>0}}\left(\frac{M_j}{N_j}\right)^{1/j}<+\infty,
\end{equation}
which is crucial in the ultradifferentiable setting. In \eqref{weightsequencrelation} instead resp. in addition of an uniform constant $A$ the geometrically growing factor $j\mapsto c^j$ is appearing with a parameter $c>0$. This difference is connected to a dilatation in the function argument in the weight $v_M$; see \eqref{powermultisequequ} and \eqref{powermultisequequ1}. It turns out, see Proposition \ref{charactprop}, that \eqref{weightsequencrelation} is becoming relevant for weighted spaces given by {\itshape non-increasing} resp. {\itshape non-decreasing systems of weights;} see e.g. \cite[Sect. 2]{BierstedtBonetTaskinen98}. On the other hand, when treating such spaces we have to replace \hyperlink{ompreceq}{$\preceq$} by a more general growth relation involving a parameter; see \eqref{bigOdilarelation}.

Moreover, in view of this observation two kinds of weight systems can be considered ''naturally'' also for abstractly given weights. We can either consider $t\mapsto v(ct)$, the \emph{dilatation-type,} motivated by the aforementioned explanations, or $t\mapsto v(t)^c$, called the \emph{exponential-type.} In both cases $c>0$ is the parameter and means, for $v_M$, that $c$ either appears ''inside'' the function argument of $\omega_M$ as a dilatation/stretching parameter resp. ''outside the function'' as a multiplicative parameter. The latter case requires another generalization of \hyperlink{ompreceq}{$\preceq$}; see \eqref{bigOexprelation}. All precise definitions can be found in Section \ref{basicweightsection}.\vspace{6pt}

In Section \ref{equsect} we are able to prove characterizations for the inclusion (resp. equality) of weighted spaces given by dilatation- and exponential-type systems. Again, this can be purely expressed in terms of given $M$ and $N$. First, we are treating systems which are based on a weight sequence and then reduce the general weight function case to the weight sequence case by involving the associated weight sequence; we refer to Theorems \ref{weightholombysequcharact}, \ref{weightholombysequcharactpower}, \ref{weightholombyfctcharact} and \ref{weightholomcharactpower}.\vspace{6pt}

Finally, based on these results in Section \ref{pointwisesection} we study and characterize closedness under the \emph{bilinear point-wise multiplication operator} $\mathfrak{m}: (f,g)\mapsto f\cdot g$ acting on such weighted classes. Again, the weight function case is reduced to the weight sequence case and in the weight sequence setting point-wise multiplication yields naturally the appearance of the so-called {\itshape convolved sequence;} see Section \ref{convolvedsect}. On the one hand, it follows immediately that the exponential-type is always closed under the action of $\mathfrak{m}$; see Proposition \ref{poinwisepowerclear}. On the other hand, in the dilatation-type weight sequence case the characterizing condition for $M$ is precisely the classical {\itshape moderate growth,} i.e. $(M.2)$ in \cite{Komatsu73}, resp. the corresponding condition \eqref{om6forv} for abstractly given $v$. We refer to the main statements Theorems \ref{convolutorthm} and \ref{convolutorthmweightfct}.

However, our characterization for weight systems is in contrast to the fact that a weighted class defined by a single weight function $v$ (see \eqref{weightedclassdef}) is not closed under $\mathfrak{m}$; see Proposition \ref{Pepecounterproof} and Corollary \ref{nonpointwiseclosed}. Moreover, note that concerning the behavior of $\mathfrak{m}$ acting on $H^{\infty}_v(G)$ there is a difference between $G=\CC$ and $G=\mathbb{D}$; in \cite{BonetDomanskiLindstroem99} the case $H^{\infty}_v(\mathbb{D})$ has been studied in detail.\vspace{6pt}

\textbf{Acknowledgements.} The author of this article thanks the anonymous referee for the careful reading and the valuable suggestions which have improved and clarified the presentation of the results. More precisely, in the detailed report the referee has remarked that the equivalence in \eqref{principalequ} follows from \cite{BonetDomanskiLindstroemTaskinen98} and so we have included this fact (main result Theorem \ref{essentialweightcharactnew}). Moreover, it has been commented that exponential-type weight systems, in particular expressed in terms of a sequence $M$, and crucial growth conditions summarized in Section \ref{modweightsection} have already been studied by Braun, Meise, and Taylor in the 1980's; see Section \ref{refereecommentsect} and the citations there.

Moreover, the author wishes to thank Prof. Jos\'{e} Bonet Solves from the Universitat Polit\`{e}cnica de Val\`{e}ncia for helpful and clarifying discussions and explanations during the preparation of this work. In particular, for forwarding and explaining (the proof of) Proposition \ref{Pepecounterproof} and the recent article \cite{Bonet2022survey}.

\section{Weights, weighted spaces and growth conditions}\label{weightsection}

\subsection{Basic convention and comments}
We use the notation $\NN=\{0,1,\dots,\}$ and $\NN_{>0}=\{1,2,\dots\}$.\vspace{6pt}

We call a weight $v$ {\itshape normalized} when $v(t)=1$ for all $t\in[0,1]$. A normalized $v$ satisfies $v(t)\le 1$ for all $t\ge 0$.

When both $v$, $w$ are (normalized) weight functions, then the product $v\cdot w$ is so, too.

For defining the class in \eqref{weightedclassdef} normalization can be assumed w.l.o.g.: Otherwise replace $v$ by a normalized weight, say $v^n$, and such that $v(t)=v^n(t)$ for all large $t>1$. This gives $H^{\infty}_v=H^{\infty}_{v^n}$ as l.c.v.s.

Since $w(t),v(t)\neq 0$ for all $t$, \eqref{bigOrelation} precisely means that
$$\exists\;C\ge 1\;\forall\;t\ge 0:\;\;\;w(t)\le Cv(t).$$

\subsection{Weight systems and corresponding weighted classes}\label{basicweightsection}
Let $\underline{\mathcal{V}}=(v_n)_{n\in\NN_{>0}}$ be a non-increasing sequence of weights, i.e. $v_n\ge v_{n+1}$ for all $n$. Then we define the (LB)-space $H^{\infty}_{\underline{\mathcal{V}}}$ of the Banach spaces $H^{\infty}_{v_n}$, so
\begin{equation}\label{indlim}
H^{\infty}_{\underline{\mathcal{V}}}:=\varinjlim_{n\rightarrow\infty}H^{\infty}_{v_n}.
\end{equation}
Analogously, if $\overline{\mathcal{V}}=(v_n)_{n\in\NN_{>0}}$ is a non-decreasing sequence of weights, i.e. $v_n\le v_{n+1}$ for all $n$, then we define the Fr\'{e}chet-space $H^{\infty}_{\overline{\mathcal{V}}}$ to be the projective limit of the Banach spaces $H^{\infty}_{v_n}$, so
\begin{equation}\label{projlim}
H^{\infty}_{\overline{\mathcal{V}}}:=\varprojlim_{n\rightarrow\infty}H^{\infty}_{v_n}.
\end{equation}
Let $v$ be a weight function and $c>0$, then set
\begin{equation}\label{parameterweights}
v_c: t\mapsto v(ct),\hspace{20pt}v^c: t\mapsto v(t)^c,
\end{equation}
and write again $v$ instead of $v_1$ resp. $v^1$. We summarize:

\begin{itemize}
\item[$(i)$] $v$ is a weight function if and only if some/each weight $v_c$ resp. if and only if some/each $v^c$ is so.

\item[$(ii)$] For any $d\ge c>0$ we have the continuous inclusions
$$H^{\infty}_{v_c}(\CC)\subseteq H^{\infty}_{v_d},\hspace{15pt}H^{\infty}_{v^c}\subseteq H^{\infty}_{v^d}.$$
Note that $v^d(t)\le v^c(t)$ for all large $t$ since eventually $v(t)\le 1$ (and for all $t\ge 0$ if $v$ is normalized).

\item[$(iii)$] $v\hyperlink{ompreceq}{\preceq}w$ if and only if $v_c\hyperlink{ompreceq}{\preceq}w_c$ if and only if $v^c\hyperlink{ompreceq}{\preceq}w^c$ for some/each $c>0$.
\end{itemize}

Consider the corresponding weight systems
\begin{equation}\label{weightsystems}
\underline{\mathcal{V}}_{\mathfrak{c}}:=(v_c)_{c\in\NN_{>0}},\hspace{20pt}\overline{\mathcal{V}}_{\mathfrak{c}}:=(v_{\frac{1}{c}})_{c\in\NN_{>0}}.
\end{equation}
Thus $\underline{\mathcal{V}}_{\mathfrak{c}}$ is a non-increasing sequence of weights, whereas $\overline{\mathcal{V}}_{\mathfrak{c}}$ is a non-decreasing sequence. Similarly, if $v\le 1$, then consider
\begin{equation}\label{weightsystems1}
\underline{\mathcal{V}}^{\mathfrak{c}}:=(v^c)_{c\in\NN_{>0}},\hspace{20pt}\overline{\mathcal{V}}^{\mathfrak{c}}:=(v^{\frac{1}{c}})_{c\in\NN_{>0}}.
\end{equation}
Since $v\le 1$ the set $\underline{\mathcal{V}}^{\mathfrak{c}}$ is a non-increasing sequence of weights, whereas $\overline{\mathcal{V}}^{\mathfrak{c}}$ is a non-decreasing sequence. The assumption $v\le 1$ holds if $v$ is normalized.

\begin{definition}
The weight systems in \eqref{weightsystems} are called to be of {\itshape dilatation-type} and the systems in \eqref{weightsystems1} of {\itshape exponential-type.}
\end{definition}

Let us write $v\hypertarget{ompreceqc}{\preceq_{\mathfrak{c}}}w$ if
\begin{equation}\label{bigOdilarelation}
\exists\;c\ge 1:\;\;\;v\hyperlink{ompreceq}{\preceq}w_c,
\end{equation}
and $v\hypertarget{simc}{\sim_{\mathfrak{c}}}w$, if
$$v\hyperlink{ompreceqc}{\preceq_{\mathfrak{c}}}w\hspace{15pt}\text{and}\hspace{15pt}w\hyperlink{ompreceqc}{\preceq_{\mathfrak{c}}}v.$$ Similarly, write $v\hypertarget{ompreceqpowc}{\preceq^{\mathfrak{c}}}w$ if
\begin{equation}\label{bigOexprelation}
\exists\;c\ge 1:\;\;\;v\hyperlink{ompreceq}{\preceq}w^c,
\end{equation}
and $v\hypertarget{simpowc}{\sim^{\mathfrak{c}}}w$, if
$$v\hyperlink{ompreceqc}{\preceq^{\mathfrak{c}}}w\hspace{15pt}\text{and}\hspace{15pt}w\hyperlink{ompreceqc}{\preceq^{\mathfrak{c}}}v.$$
Note:
\begin{itemize}
\item[$(*)$] Both $v\hyperlink{ompreceqc}{\preceq_{\mathfrak{c}}}w$ and $v\hyperlink{ompreceqpowc}{\preceq^{\mathfrak{c}}}w$ are weaker than $v\hyperlink{ompreceq}{\preceq}w$.

\item[$(*)$] $v\hyperlink{ompreceqc}{\preceq_{\mathfrak{c}}}w$ resp. $v\hyperlink{ompreceqc}{\preceq^{\mathfrak{c}}}w$ if and only if $$\exists\;c\ge 1\;\forall\;a>0:\;\;\;v_a\hyperlink{ompreceq}{\preceq}w_{ac}\;\;\;\text{resp.}\;\;\;v^a\hyperlink{ompreceq}{\preceq}w^{ac}.$$
\end{itemize}

These observations imply the following:

\begin{proposition}\label{firstweighfctprop}
Let $v$ and $w$ be (normalized) weight functions. If $v\hyperlink{ompreceqc}{\preceq_{\mathfrak{c}}}w$, then
$$H^{\infty}_{\underline{\mathcal{V}}_{\mathfrak{c}}}\subseteq H^{\infty}_{\underline{\mathcal{W}}_{\mathfrak{c}}},\hspace{20pt}H^{\infty}_{\overline{\mathcal{V}}_{\mathfrak{c}}}\subseteq H^{\infty}_{\overline{\mathcal{W}}_{\mathfrak{c}}},$$
and if $v\hyperlink{ompreceqpowc}{\preceq^{\mathfrak{c}}}w$, then
$$H^{\infty}_{\underline{\mathcal{V}}^{\mathfrak{c}}}\subseteq H^{\infty}_{\underline{\mathcal{W}}^{\mathfrak{c}}},\hspace{20pt}H^{\infty}_{\overline{\mathcal{V}}^{\mathfrak{c}}}\subseteq H^{\infty}_{\overline{\mathcal{W}}^{\mathfrak{c}}},$$
with continuous inclusions.
\end{proposition}

We close with the following observation:

\begin{remark}\label{strongpowerremark}
\emph{Let $v$ be a (normalized) weight function, then
\begin{equation}\label{strongweightsrelation1}
\forall\;d>c>0:\;\;\;\lim_{t\rightarrow+\infty}\frac{v^c(t)}{v^d(t)}=\lim_{t\rightarrow+\infty}v^{c-d}(t)=+\infty.
\end{equation}
However, in general $\lim_{t\rightarrow+\infty}\frac{v_c(t)}{v_d(t)}=+\infty$ is not clear.}
\end{remark}

\subsection{Weight sequences and conditions}
Let $M=(M_j)_j\in\RR_{>0}^{\NN}$. The corresponding sequence of quotients is denoted by $\mu_j:=\frac{M_j}{M_{j-1}}$, $j\ge 1$, and set $\mu_0:=1$. $M$ is called {\itshape normalized} if $1=M_0\le M_1$ and $M$ is called {\itshape log-convex} if
$$\forall\;j\in\NN_{>0}:\;M_j^2\le M_{j-1} M_{j+1},$$
equivalently if the sequence $(\mu_j)_j$ is non-decreasing. If $M_0=1$ and $M$ is log-convex, then $(M_j)^{1/j}\le\mu_j$ for all $j\in\NN_{>0}$.

\begin{definition}\label{defweightsequ}
$M$ is called a {\itshape weight sequence} if
$$1=M_0\;\;\;\text{and}\;\;\;\lim_{j\rightarrow+\infty}(M_j)^{1/j}=+\infty.$$
\end{definition}

Moreover consider the set
$$\hypertarget{LCset}{\mathcal{LC}}:=\{M\in\RR_{>0}^{\NN}:\;M\;\text{is normalized, log-convex},\;\lim_{j\rightarrow+\infty}(M_j)^{1/j}=+\infty\},$$
i.e. $M\in\hyperlink{LCset}{\mathcal{LC}}$ if and only if $M$ is a log-convex weight sequence with $M_1\ge 1$. Note that each $M\in\hyperlink{LCset}{\mathcal{LC}}$ is non-decreasing since $\mu_j\ge 1$ for all $j$ and $j\mapsto\mu_j$ is non-decreasing.\vspace{6pt}

$M$ (with $M_0=1$) has condition {\itshape moderate growth}, denoted by \hypertarget{mg}{$(\text{mg})$}, if
$$\exists\;C\ge 1\;\forall\;j,k\in\NN:\;M_{j+k}\le C^{j+k} M_j M_k.$$
In \cite{Komatsu73} this is denoted by $(M.2)$ and also known under the name {\itshape stability under ultradifferential operators.}\vspace{6pt}

Let $M,N\in\RR_{>0}^{\NN}$, we write $M\le N$ if $M_j\le N_j$ for all $j$ and $M\hypertarget{preceq}{\preceq}N$ if \eqref{weightsequencrelation}, i.e.
$$\sup_{j\in\NN_{>0}}\left(\frac{M_j}{N_j}\right)^{1/j}<+\infty.$$

$M$ and $N$ are called {\itshape equivalent}, written $M\hypertarget{approx}{\approx}N$, if
$$M\hyperlink{preceq}{\preceq}N\hspace{15pt}\text{and}\hspace{15pt}N\hyperlink{preceq}{\preceq}M.$$
Property \hyperlink{mg}{$(\on{mg})$} is clearly preserved under \hyperlink{approx}{$\approx$}.

\begin{remark}\label{LCapproxrem}
\emph{Let $M$ be a log-convex weight sequence then there exists an equivalent $N\in\hyperlink{LCset}{\mathcal{LC}}$ such that $\mu_j=\nu_j$ for all $j$ sufficiently large: Since $\lim_{j\rightarrow+\infty}\mu_j=+\infty$ and $(\mu_j)_j$ is non-decreasing, there exists $j_0\in\NN_{>0}$ such that $\mu_j\ge\mu_{j_0}>1$ for all $j\ge j_0$. Then set}
$$\nu_j:=1,\;\;\;j<j_0,\hspace{15pt}\text{and}\hspace{15pt}\nu_j:=\mu_j,\;\;\;j\ge j_0,$$
\emph{and consequently}
\begin{equation}\label{LCapproxremequ}
\exists\;C\ge 1\;\forall\;j\in\NN:\;\;\;\frac{1}{C}N_j\le M_j\le CN_j.
\end{equation}
\end{remark}

\subsection{Associated weight function}\label{assofctsect}
Let $M\in\RR_{>0}^{\NN}$ (with $M_0=1$), then the {\itshape associated weight function} $\omega_M: \RR_{\ge 0}\rightarrow\RR\cup\{+\infty\}$ is defined by
\begin{equation*}\label{assofunc}
\omega_M(t):=\sup_{j\in\NN}\log\left(\frac{t^j}{M_j}\right)\;\;\;\text{for}\;t\in\RR_{>0},\hspace{30pt}\omega_M(0):=0.
\end{equation*}
For an abstract introduction of the associated function we refer to \cite[Chapitre I]{mandelbrojtbook}, see also \cite[Definition 3.1]{Komatsu73} and the more recent work \cite{regularnew}.

If $\liminf_{j\rightarrow+\infty}(M_j)^{1/j}>0$, then $\omega_M(t)=0$ for sufficiently small $t$, since $\log\left(\frac{t^j}{M_j}\right)<0\Leftrightarrow t<(M_j)^{1/j}$ holds for all $j\in\NN_{>0}$. (In particular, if $M_j\ge 1$ for all $j\in\NN$, then $\omega_M$ is vanishing on $[0,1]$.) Moreover, under this assumption $t\mapsto\omega_M(t)$ is a continuous non-decreasing function, which is convex in the variable $\log(t)$ and tends faster to infinity than any $\log(t^j)$, $j\ge 1$, as $t\rightarrow+\infty$. $\lim_{j\rightarrow+\infty}(M_j)^{1/j}=+\infty$ implies that $\omega_M(t)<+\infty$ for each finite $t$ which shall be considered as a basic assumption for defining $\omega_M$. In particular, this holds true provided $M$ is a weight sequence.\vspace{6pt}

If $M$ is a log-convex weight sequence, then let us introduce the {\itshape counting function}
\begin{equation}\label{counting}
\Sigma_{M}(t):=|\{j\ge 1:\;\;\;\mu_j\le t\}|.
\end{equation}
In this case we can compute $M$ by involving $\omega_M$ as follows, see \cite[Chapitre I, 1.4, 1.8]{mandelbrojtbook} and also \cite[Prop. 3.2]{Komatsu73}:
\begin{equation}\label{Prop32Komatsu}
M_j=\sup_{t\ge 0}\frac{t^j}{\exp(\omega_{M}(t))},\;\;\;j\in\NN.
\end{equation}
If $M$ is a weight sequence which is not log-convex, then the right-hand side of \eqref{Prop32Komatsu} yields $M^{\on{lc}}_j$ with $M^{\on{lc}}$ denoting the {\itshape log-convex minorant} of $M$ and $M^{\on{lc}}$ is a log-convex weight sequence.

Moreover, one has
\begin{equation*}\label{assovanishing}
\forall\;t\in[0,\mu_1]:\;\;\;\omega_M(t)=0,
\end{equation*}
which follows by the known integral representation formula (see \cite[1.8. III]{mandelbrojtbook} and also \cite[$(3.11)$]{Komatsu73})
\begin{equation}\label{assointrepr}
\omega_M(t)=\int_0^t\frac{\Sigma_M(u)}{u}du=\int_{\mu_1}^t\frac{\Sigma_M(u)}{u}du.
\end{equation}
In particular, if $M\in\hyperlink{LCset}{\mathcal{LC}}$, then $\omega_M$ vanishes on $[0,1]$.

Finally, for any weight sequence $M$ we have
\begin{equation}\label{lcminorantnodifference}
\omega_{M^{\on{lc}}}=\omega_{M}.
\end{equation}

\subsection{Weighted spaces of entire functions defined in terms of $\omega_M$}
For any weight sequence $M$ (in the sense of Definition \ref{defweightsequ}) and for any parameter $c>0$ we set
\begin{equation}\label{weights}
v_{M,c}(t):=\exp(-\omega_M(ct)),\;\;\;v^c_M(t):=\exp(-c\omega_M(t)),\;\;\;t\ge 0.
\end{equation}
Note that each $v_{M,c}$ and $v^c_M$ is a weight by the properties of the associated weight function. Then define
\begin{align*}
H^{\infty}_{v_{M,c}}&:=\{f\in H(\CC): \|f\|_{v_{M,c}}:=\sup_{z\in\CC}|f(z)|v_{M,c}(|z|)<+\infty\}
\\&
=\{f\in H(\CC): \|f\|_{v_{M,c}}:=\sup_{z\in\CC}|f(z)|\exp(-\omega_M(c|z|))<+\infty\},
\end{align*}
and the class $H^{\infty}_{v^c_M}$ is introduced accordingly. Recall that $M\in\hyperlink{LCset}{\mathcal{LC}}$ has been the standard assumption in \cite{solidassociatedweight}.

\begin{remark}\label{LCwlogrem}
\emph{We comment on the requirements for $M$:}
\begin{itemize}
\item[$(*)$] \emph{If $M$ is a weight sequence, then in view of \eqref{lcminorantnodifference} there is no difference between the weighted classes defined in terms of $M$ or $M^{\on{lc}}$ and thus one can assume w.l.o.g. that the defining sequence is log-convex.}

\item[$(*)$] \emph{Let $M$ be a log-convex weight sequence. By Remark \ref{LCapproxrem} there exists $N\in\hyperlink{LCset}{\mathcal{LC}}$ such that \eqref{LCapproxremequ} holds true and so by definition
$$\exists\;D\ge 0\;\forall\;t\ge 0:\;\;\;\omega_M(t)-D\le\omega_N(t)\le\omega_M(t)+D,$$
which yields as l.c.v.s.
$$\forall\;c>0:\;\;\;H^{\infty}_{v_{M,c}}=H^{\infty}_{v_{N,c}},\hspace{15pt}H^{\infty}_{v^c_M}=H^{\infty}_{v^c_N}.$$}
\end{itemize}
\end{remark}

\begin{remark}\label{LCwlogrem1}
\emph{We summarize more consequences for weight sequences $M$ and $N$:}

\begin{itemize}
\item[$(i)$] \emph{For $c=1$ we write $v_M$ instead of $v_{M,1}$ resp. $v^1_M$. $v_M$ is normalized if $\omega_M$ vanishes on $[0,1]$, in particular this holds for any $M\in\hyperlink{LCset}{\mathcal{LC}}$.}

\item[$(ii)$] \emph{One has that
$$\forall\;d\ge c>0:\;\;\;H^{\infty}_{v_{M,c}}\subseteq H^{\infty}_{v_{M,d}},\hspace{15pt}H^{\infty}_{v^c_M}\subseteq H^{\infty}_{v^d_M},$$
with continuous inclusions.}

\item[$(iii)$] \emph{For the weight systems in \eqref{weightsystems} and \eqref{weightsystems1} we set ($M$ is fixed)}
\begin{equation}\label{weightsequweightsystem}
\underline{\mathcal{M}}_{\mathfrak{c}}:=(v_{M,c})_{c\in\NN_{>0}},\hspace{10pt}\overline{\mathcal{M}}_{\mathfrak{c}}:=(v_{M,\frac{1}{c}})_{c\in\NN_{>0}},\hspace{10pt}\underline{\mathcal{M}}^{\mathfrak{c}}:=(v^c_M)_{c\in\NN_{>0}},\hspace{10pt}\overline{\mathcal{M}}^{\mathfrak{c}}:=(v^{\frac{1}{c}}_{M})_{c\in\NN_{>0}},
\end{equation}
\emph{and the (LB)-spaces from \eqref{indlim} resp. the Fr\'{e}chet-spaces from \eqref{projlim} are defined accordingly.}

\item[$(iv)$] \emph{By Remark \ref{strongpowerremark} one has}
$$\forall\;d>c>0:\;\;\;\lim_{t\rightarrow+\infty}\frac{v^c_M(t)}{v^d_M(t)}=+\infty.$$

\item[$(v)$] \emph{The following are equivalent by definition:}
\begin{itemize}
\item[$(*)$] \emph{$\omega_M(t)=O(\omega_N(t))$, i.e. $\omega_N\hyperlink{ompreceq}{\preceq}\omega_M$,}

\item[$(*)$] \emph{for some $c\ge 1$ we get $v_M\hyperlink{ompreceq}{\preceq}v_N^c$, i.e. $v_M\hyperlink{ompreceqpowc}{\preceq^{\mathfrak{c}}}v_N$.}
\end{itemize}
\end{itemize}
\end{remark}

We prove now a first characterization in terms of \eqref{strongNMrelation}.

\begin{proposition}\label{strongNMrelationlemma}
Let $M$ and $N$ be weight sequences and assume that $N$ is log-convex. Then the following are equivalent:
\begin{itemize}
\item[$(i)$] The sequences satisfy \eqref{strongNMrelation}, i.e.
$$\exists\;A\ge 1\;\forall\;j\in\NN:\;\;\;N_j\le AM_j.$$

\item[$(ii)$] The weights satisfy
$$v_M\hyperlink{ompreceq}{\preceq}v_N.$$
\end{itemize}
\end{proposition}
For $(i)\Rightarrow(ii)$ also log-convexity for $N$ is not required necessarily and recall that $v_M\hyperlink{ompreceq}{\preceq}v_N$ implies $H^{\infty}_{v_M}\subseteq H^{\infty}_{v_N}$ (with continuous inclusion).

\demo{Proof}
$(i)\Rightarrow(ii)$ By assumption and definition $\omega_M(t)\le\omega_N(t)+\log(A)$ for all $t\ge 0$ follows and so $(ii)$ is immediate.

$(ii)\Rightarrow(i)$ By assumption there exists some $C\ge 1$ such that $\exp(-\omega_N(t))=v_N(t)\le Cv_M(t)=C\exp(-\omega_M(t))$ for all $t\ge 0$. Then \eqref{Prop32Komatsu} applied to both weights yields $N_j\le CM^{\on{lc}}_j\le CM_j$ for all $j\in\NN$, i.e. \eqref{strongNMrelation} with $A=C$.
\qed\enddemo

The next statement should be compared with Remark \ref{strongpowerremark}.

\begin{lemma}\label{strongdilaremark}
Let $M$ be a log-convex weight sequence. Then
\begin{equation}\label{strongweightsrelation}
\forall\;d>c>0:\;\;\;\lim_{t\rightarrow+\infty}\frac{v_{M,c}(t)}{v_{M,d}(t)}=+\infty.
\end{equation}
\end{lemma}

\demo{Proof}
Let $d>c>0$, then by the definition of $\Sigma_M$ in \eqref{counting} we get
$$\omega_M(dt)-\omega_M(ct)=\int_{ct}^{dt}\frac{\Sigma_M(u)}{u}du\ge\Sigma_M(ct)\int_{ct}^{dt}\frac{1}{u}du=\Sigma_M(ct)\log(d/c)\rightarrow+\infty.$$
\qed\enddemo

The appearing parameter $c>0$ in \eqref{weights} in the dilatation-type structure is closely connected to and naturally appearing in the weight sequence approach.

\begin{remark}\label{powermultisequ}
\emph{Let $M$ be a given weight sequence.}

\begin{itemize}
\item[$(a)$] \emph{For any $c>0$ we set $M^c:=(c^jM_j)_{j\in\NN}$ with corresponding quotient sequence $\mu^c_j=c\mu_j$ for $j\ge 1$ (and again $\mu^c_0:=1$). Then $1=M^c_0$ and $\lim_{j\rightarrow+\infty}(M^c_j)^{1/j}=+\infty$, i.e. $M^c$ is a weight sequence as well and $M$ is log-convex if and only if some/each $M^c$ is so. We also have}
    \begin{equation}\label{powermultisequequ}
    \forall\;c>0\;\forall\;t\ge 0:\;\;\;M^c\hyperlink{approx}{\approx}M,\hspace{15pt}\omega_{M^c}(t)=\omega_M\left(\frac{t}{c}\right),
    \end{equation}
    \emph{which gives}
\begin{equation}\label{powermultisequequ1}
\forall\;c>0\;\forall\;t\ge 0:\;\;\;v_{M,c}(t)=v_{M^{\frac{1}{c}}}(t).
\end{equation}
\emph{If $c<1$, then even if $M$ is normalized in general normalization is not clear for $M^c$.}

\item[$(b)$] \emph{For arbitrary $c>0$ the spaces $H^{\infty}_{v_{M,c}}$ and $H^{\infty}_{v_{M}}$ are isometrically isomorphic, written $H^{\infty}_{v_{M,c}}\cong H^{\infty}_{v_{M}}$: This follows by taking into account the mappings $f(z)\mapsto f(zc)$ and $f(z)\mapsto f(z/c)$, see also \cite[Cor. 4.5, Rem. 7.1]{solidassociatedweight}. Analogously, this comment holds true for an abstractly given weight $v$.}
\end{itemize}
\end{remark}

Let $f_k(z):=z^k$, then for all $k\in\NN$ and any weight sequence $M$ we have $f_k\in H^{\infty}_{\overline{\mathcal{M}}_{\mathfrak{c}}}, H^{\infty}_{\overline{\mathcal{M}}^{\mathfrak{c}}}$, because $|f_k(z)|v_{M,c}(|z|)=|z|^kv_M(c|z|)\rightarrow 0$ and $|f_k(z)|v^c_M(|z|)=|z|^kv^c_M(|z|)\rightarrow 0$ as $|z|\rightarrow+\infty$ for all $c>0$. For this recall that $v_M$ is rapidly decreasing. The same statement holds true if $v_M$ is replaced by any abstractly given weight $v$: In fact $H_v^{\infty}$ contains all polynomials if and only if $v$ is rapidly decreasing; see \cite[p. 3]{Bonet2022survey}.

\subsection{On the class of (log-)convex weights $v$}\label{convexsection}
We recall some statements from \cite[Sect. 6]{solidassociatedweight}. Let $v:[0,+\infty)\rightarrow(0,+\infty)$ be a normalized (radial) weight function. Then set
\begin{equation}\label{omegafromv}
\omega^v(t):=-\log(v(t)),\;\;\;t\in[0,+\infty),
\end{equation}
and we get (cf. \cite[Lemma 6.1]{solidassociatedweight}):
\begin{itemize}
\item[$(*)$] $\omega^v:[0,+\infty)\rightarrow[0,+\infty)$ is continuous and non-decreasing,

\item[$(*)$] $\lim_{t\rightarrow+\infty}\omega^v(t)=+\infty$ and

\item[$(*)$] $\omega^v(t)=0$ for $t\in[0,1]$ (normalization).

\item[$(*)$] The fact that $v$ is rapidly decreasing is equivalent to
$$\hypertarget{om3}{(\omega_3)}:\;\;\;\log(t)=o(\omega^v(t))\hspace{15pt}t\rightarrow+\infty.$$

\item[$(*)$] We have that
\begin{equation}\label{vconvexity}
t\mapsto\omega^v(e^t)(=-\log(v(e^t)))\;\;\;\text{is convex on}\;\RR,
\end{equation}
if and only if $\varphi_{\omega^v}: t\mapsto\omega^v(e^t)$ is convex.
\end{itemize}

\hyperlink{om3}{$(\omega_3)$} and \eqref{vconvexity} are standard assumptions in the theory of ultradifferentiable functions defined by so-called {\itshape Braun-Meise-Taylor weight functions} $\omega$; see \cite{BraunMeiseTaylor90}. \hyperlink{om3}{$(\omega_3)$} is named after \cite{dissertation} and there \eqref{vconvexity} is abbreviated with $(\omega_4)$. Note that $v$ has \eqref{vconvexity} if and only if $v_c$ resp. $v^c$ for some/each $c>0$ does so and $\omega^{v_c}(t)=-\log(v_c(t))=-\log(v(tc))$, whereas $\omega^{v^c}(t)=-\log(v^c(t))=-c\log(v(t))=-c\omega^v(t)$.\vspace{6pt}

Conversely, we write that $\omega: [0,+\infty)\rightarrow[0,+\infty)$ has $\hypertarget{om0}{(\omega_0)}$ if $\omega$ is continuous, non-decreasing, $\omega(t)=0$ for all $t\in[0,1]$ (normalization) and $\lim_{t\rightarrow+\infty}\omega(t)=+\infty$. We put
$$\hypertarget{omset0}{\mathcal{W}_0}:=\{\omega:[0,\infty)\rightarrow[0,\infty): \omega\;\text{has}\;\hyperlink{om0}{(\omega_0)},\hyperlink{om3}{(\omega_3)},\eqref{vconvexity}\},$$
and for any $\omega\in\hyperlink{omset0}{\mathcal{W}_0}$ we set
\begin{equation}\label{vfromomega}
v^{\omega}(t):=\exp(-\omega(t)).
\end{equation}
Then $v^{\omega}$ is a normalized weight function such that \eqref{vconvexity} holds true. Summarizing, we have shown:

\begin{lemma}\label{vversusomega}
The set of normalized weight functions $v$ having \eqref{vconvexity} is isomorphic to the set \hyperlink{omset0}{$\mathcal{W}_0$} via the mappings
$$v\mapsto\omega^v,\hspace{20pt}\omega\mapsto v^{\omega}.$$
\end{lemma}

In particular, Lemma \ref{vversusomega} can be applied to $\omega\equiv\omega_M$ provided $M\in\hyperlink{LCset}{\mathcal{LC}}$ since in this case $\omega_M\in\hyperlink{omset0}{\mathcal{W}_0}$; see e.g. \cite[Lemma 3.1]{sectorialextensions1} and the citations there. In this situation $v^{\omega_M}\equiv v_M$ and $\omega^{v_M}\equiv\omega_M$.

\begin{definition}\label{admissdef}
We call a weight $v$ {\itshape convex} if it satisfies \eqref{vconvexity}.
\end{definition}

Note that in \cite[p. 5, $(4)$]{Bonet2022survey} a weight with property \eqref{vconvexity} has been called {\itshape log-convex;} this has also been done in \cite{AbakumovDoubtsov15} and \cite{AbakumovDoubtsov18}.

\subsection{From weight functions to weight sequences}\label{weighttosequsect}
The following is also taken from \cite[Sect. 6]{solidassociatedweight}. Let $v$ be a normalized and convex weight and let $\omega^v$ be given by \eqref{omegafromv}. In \cite[$(6.3)$]{solidassociatedweight} we have introduced and studied the sequence $M^v=(M^v_j)_{j\in\NN}$ defined by
\begin{equation}\label{vBMTweight1equ1}
M^v_j:=\sup_{t>0}\frac{t^j}{\exp(\omega^v(t))}=\sup_{t>0}t^jv(t),\;\;\;j\in\NN.
\end{equation}
Note that this expression has already been introduced on \cite[p. 157]{BierstedtBonetTaskinen98}.

We are also interested in the weight $t\mapsto\exp(-\omega_{M^v}(t))(=v_{M^v}(t))$. This notation is suggested by \eqref{weights} but confusing. Thus when explicitly using this function we denote the abstractly given weight with a different symbol than $v$, say $w$ or $u$, e.g. $v_{M^u}(t)=\exp(-\omega_{M^u}(t))$.

When $M\equiv M^u$, then in \eqref{weightsequweightsystem} replace $\mathcal{M}$ by $\mathcal{M}^u$.\vspace{6pt}

We gather several properties; see \cite[Prop. 6.2]{solidassociatedweight}:

\begin{itemize}
\item[$(i)$] $M^u\in\hyperlink{LCset}{\mathcal{LC}}$,

\item[$(ii)$] $\omega_{M^u}\hyperlink{sim}{\sim}\omega^u$, more precisely $\exists\;A\ge 1\;\forall\;t\ge 0:$
\begin{equation}\label{omegavequiv}
\frac{1}{A}(v_{M^{u}}(t))^{2}=\frac{1}{A}\exp(-2\omega_{M^u}(t))\le u(t)=\exp(-\omega^u(t))\le\exp(-\omega_{M^u}(t))=v_{M^u}(t),
\end{equation}
consequently
\begin{equation}\label{omegavequivnew}
\exists\;A\ge 1\;\forall\;c>0\;\forall\;t\ge 0:\;\;\;\frac{1}{A}v^2_{M^u,c}(t)\le u_c(t)\le v_{M^u,c}(t),\;\;\;\frac{1}{A}v^{2c}_{M^u}(t)\le u^c(t)\le v^c_{M^u}(t),
\end{equation}
which implies
\begin{equation}\label{omegavequivnewnew}
\forall\;c>0:\;\;\;H^{\infty}_{v_{M^u},c}\subseteq H^{\infty}_{u,c}\subseteq H^{\infty}_{v^2_{M^u},c},\hspace{15pt}H^{\infty}_{v^c_{M^u}}\subseteq H^{\infty}_{u^c}\subseteq H^{\infty}_{v^{2c}_{M^u}},
\end{equation}
with continuous inclusions. Recall that the second part in \eqref{omegavequivnew} precisely means $v_{M^u}\hyperlink{simpowc}{\sim^{\mathfrak{c}}}u$.

\item[$(iii)$] If $u\equiv v_M$ with $M\in\hyperlink{LCset}{\mathcal{LC}}$, and so $\omega^u\equiv\omega_M$, then $M^u\equiv M$. Moreover, for any $c>0$ by definition $(M^u)^{\frac{1}{c}}\equiv(c^jM^u_j)_{j\in\NN}$.
\end{itemize}

We comment on the consequences of the standard requirements for a given weight.

\begin{remark}\label{convnoneed}
\emph{A careful inspection of the proof of \cite[Prop. 6.2]{solidassociatedweight}, and of the used citations \cite[Thm. 4.0.3]{dissertation}, \cite[Lemma 5.7]{compositionpaper}, gives that convexity for $u$ is in the whole approach only required in order to get the first estimate in \eqref{omegavequiv} resp. \eqref{omegavequivnew}; i.e. the second inclusions in \eqref{omegavequivnewnew}. Crucially we remark that}
\begin{equation*}\label{convnoneedequ}
\omega_{M^u}\hyperlink{sim}{\sim}\omega^u\nRightarrow v_{M^u}\hyperlink{sim}{\sim}u,
\end{equation*}
\emph{which is due to the appearance of the quadratic power in the first estimate in \eqref{omegavequiv}. This fact requires a technical condition studied in Section \ref{modweightsection}.}
\end{remark}

\begin{remark}\label{assoweightsequ}
\emph{If $v$ is a function having all properties to be a normalized weight except that $v$ is rapidly decreasing, then \eqref{vBMTweight1equ1} implies that $M^v_j=+\infty$ for all sufficiently large $j$ and so the technique of associating a weight sequence is not available anymore.}\vspace{6pt}

\emph{Next let us see that $\sup_{j\in\NN_{>0}}(M^v)^{1/j}<+\infty$ if and only if $v(t)=0$ for all $t$ sufficiently large (i.e. $v$ is violating strict positivity):}\vspace{6pt}

\emph{If the supremum is finite, then by \eqref{vBMTweight1equ1} there exists $D\ge 1$ such that for all $t\ge 0$ and $j\in\NN_{>0}$ we get $t^jv(t)\le D^j\Leftrightarrow v(t)\le(D/t)^j$. Thus $v(t)=0$ for all $t>D$ follows.}

\emph{Conversely, if there exists $t_0$ with $v(t)=0$ for all $t\ge t_0$, then $M^v_j=\sup_{t>0}t^jv(t)=\max_{0\le t\le t_0}t^jv(t)\le v(0)t_0^j=t_0^j$ for all $j\in\NN$ since $v$ is non-decreasing and $v(0)=1$.}

\emph{Note that for weights $v$ such that $v(t)=0$ for all sufficiently large $t$ we get that $H_v^{\infty}=H(\CC)$ (as sets!).}
\end{remark}

The second part of \eqref{omegavequivnewnew} immediately gives:

\begin{proposition}\label{omegavequivnewnewnew}
Let $u$ be a normalized and convex weight function, $M^u$ the associated weight sequence and consider the corresponding weight systems $\underline{\mathcal{U}}^{\mathfrak{c}}$ resp. $\overline{\mathcal{U}}^{\mathfrak{c}}$ (see \eqref{weightsystems1}) and $\underline{\mathcal{M}^u}^{\mathfrak{c}}$ resp. $\overline{\mathcal{M}^u}^{\mathfrak{c}}$ (see \eqref{weightsequweightsystem}). Then as locally convex vector spaces we get
$$H^{\infty}_{\underline{\mathcal{M}^u}^{\mathfrak{c}}}= H^{\infty}_{\underline{\mathcal{U}}^{\mathfrak{c}}},\hspace{15pt}H^{\infty}_{\overline{\mathcal{M}^u}^{\mathfrak{c}}}=H^{\infty}_{\overline{\mathcal{U}}^{\mathfrak{c}}}.$$
If $u$ is not convex, then in view of Remark \ref{convnoneed} we only get the (continuous) inclusions $\subseteq$.
\end{proposition}

We close with the following important observation concerning the norm of the monomials $f_j(z)=z^j$, $j\in\NN$. For a given weight $v$ we set (see \cite[p. 5]{Bonet2022survey}):
\begin{equation}\label{Pvformula}
P_v(t):=\sup_{j\in\NN}\frac{t^j}{\|f_j\|_v},\;\;\;t\ge 0.
\end{equation}

\begin{lemma}\label{monomialnormlemma}
Let $u$ be a given normalized weight. Then
$$\forall\;j\in\NN:\;\;\;\|f_j\|_u=M^u_j,$$
and
\begin{equation}\label{monomialnormlemmaequ}
\forall\;t\ge 0:\;\;\;P_u(t)=\exp(\omega_{M^u}(t))=\frac{1}{v_{M^u}(t)}\le\frac{1}{u(t)}.
\end{equation}
If $u=v_M$ for a given log-convex weight sequence $M$, then
$$\forall\;j\in\NN:\;\;\;\|f_j\|_{v_M}=M_j,\hspace{30pt}\forall\;t\ge 0:\;\;\;P_{v_M}(t)=\exp(\omega_M(t))=\frac{1}{v_M(t)}.$$
\end{lemma}

The inequality in \eqref{monomialnormlemmaequ} has already been observed in \cite[Sect. 2.1]{AbakumovDoubtsov18}. This follows immediately by definition and also by the second half in \eqref{omegavequiv}.

\demo{Proof}
By definition and \eqref{vBMTweight1equ1} we have for all $j\in\NN$ that
$$\|f_j\|_u=\sup_{z\in\CC}|z^j|u(|z|)=\sup_{t\ge 0}t^ju(t)=M^u_j,$$
and so $P_u(t)=\sup_{j\in\NN}\frac{t^j}{\|f_j\|_u}=\sup_{j\in\NN}\frac{t^j}{M^u_j}=\exp(\omega_{M^u}(t))=\frac{1}{v_{M^u}(t)}$ is immediate.

In view of \eqref{Prop32Komatsu} the case $u=v_M$ is then clear. Note that, if $M$ is not log-convex, then we get the equality $\|f_j\|_{v_M}=M^{\on{lc}}_j$ but by Remark \ref{LCwlogrem} log-convexity is not restricting.
\qed\enddemo

\subsection{Technical growth conditions on weight functions}\label{modweightsection}
We revisit two (natural) growth conditions which have already appeared in the literature for both the ultradifferentiable and the weighted entire setting; see Section \ref{refereecommentsect} below for more comments and also Section \ref{convolvedsect}.

From $\omega_{M^u}\hyperlink{sim}{\sim}\omega^u$ it does not follow that $u\hyperlink{sim}{\sim}v_{M^u}$ is valid; recall Remark \ref{convnoneed}. In order to get more information we recall the following property:

$$\hypertarget{om6}{(\omega_6)}:\;\;\; \exists\;H\ge 1\;\forall\;t\ge 0:\;\;\;2\omega(t)\le\omega(H t)+H.$$

This condition is named after \cite{compositionpaper} and \cite{dissertation}. For abstractly given weight functions $\omega$ (in the sense of Braun-Meise-Taylor) \hyperlink{om6}{$(\omega_6)$} has appeared in \cite{BonetMeiseMelikhov07} and for a log-convex weight sequence $M$ the function $\omega_M$ has \hyperlink{om6}{$(\omega_6)$} if and only if $M$ has \hyperlink{mg}{$(\on{mg})$}; see \cite[Prop. 3.6]{Komatsu73}.

\begin{lemma}\label{om6rem}
Let $u$ be normalized and convex, then the following are equivalent:
\begin{itemize}
\item[$(a)$] $u$ satisfies
\begin{equation}\label{om6forv}
\exists\;H\ge 1\;\forall\;t\ge 0:\;\;\;u_H(t)=u(Ht)\le e^Hu^2(t).
\end{equation}

\item[$(b)$] The function $\omega^u$ (see \eqref{omegafromv}) satisfies \hyperlink{om6}{$(\omega_6)$}.

\item[$(c)$] The associated weight function $\omega_{M^u}$ satisfies \hyperlink{om6}{$(\omega_6)$}.

\item[$(d)$] The sequence $M^u$ satisfies \hyperlink{mg}{$(\on{mg})$}.
\end{itemize}
\end{lemma}

\demo{Proof}
$(a)\Leftrightarrow(b)$ follows with the same $H$ directly by definition (see \eqref{omegafromv}); $(b)\Leftrightarrow(c)$ holds by $\omega_{M^u}\hyperlink{sim}{\sim}\omega^u$ (for which convexity is required!) and since \hyperlink{om6}{$(\omega_6)$} is preserved under \hyperlink{sim}{$\sim$}; $(c)\Leftrightarrow(d)$ holds by \cite[Prop. 3.6]{Komatsu73}, see also $(iii)$ in Remark \ref{convolvesequrem}.
\qed\enddemo

{\itshape Note:}

\begin{itemize}
\item[$(*)$] \eqref{om6forv} holds if and only if there exist $H_1,H_2\ge 1$ such that $u_{H_1}(t)\le e^{H_2}u^2(t)$ for all $t\ge 0$ (take $H:=\max\{H_1,H_2\}$).

\item[$(*)$] \eqref{om6forv} holds for $u$ if and only if for some/each $u_c$ even with the same choice for $H$ and it is clearly preserved under relation \hyperlink{sim}{$\sim$}.

\item[$(*)$] The choice $H>1$ is required since for $0<H\le 1$ one has $u(t)\le u(Ht)\le e^Hu^2(t)$ and so $1\le e^Hu(t)$ for all $t\ge 0$, a contradiction.
\end{itemize}

\begin{definition}\label{admissdef1}
A weight $v$ is called to be of {\itshape moderate growth} if $v$ satisfies \eqref{om6forv}.
\end{definition}

Combining the first part of \eqref{omegavequivnewnew} and \eqref{om6forv} we get:

\begin{proposition}\label{omegavequivnewnewnewmg}
Let $u$ be a normalized and convex weight function which is of moderate growth. Let $M^u$ be the associated weight sequence and consider the corresponding weight systems
$\underline{\mathcal{U}}_{\mathfrak{c}}$ resp. $\overline{\mathcal{U}}_{\mathfrak{c}}$ (see \eqref{weightsystems}) and $\underline{\mathcal{M}^u}_{\mathfrak{c}}$  resp. $\overline{\mathcal{M}^u}_{\mathfrak{c}}$ (see \eqref{weightsequweightsystem}). Then as locally convex vector spaces we get
$$H^{\infty}_{\underline{\mathcal{U}}_{\mathfrak{c}}}=H^{\infty}_{\underline{\mathcal{M}^u}_{\mathfrak{c}}},\hspace{15pt}H^{\infty}_{\overline{\mathcal{U}}_{\mathfrak{c}}}=H^{\infty}_{\overline{\mathcal{M}^u}_{\mathfrak{c}}}.$$
\end{proposition}

For weight functions $\omega$ in the sense of Braun-Meise-Taylor also the following assumption is standard; see \cite{BraunMeiseTaylor90}:
$$\hypertarget{om1}{(\omega_1)}:\;\;\; \exists\;L\ge 1\;\forall\;t\ge 0:\;\;\;\omega(2t)\le L(\omega(t)+1).$$
This condition is named again after \cite{compositionpaper} and \cite{dissertation}. Recently \hyperlink{om1}{$(\omega_1)$} for $\omega_M$ has been characterized in terms of $M$ in \cite[Thm. 3.1]{subaddlike}: If $M\in\hyperlink{LCset}{\mathcal{LC}}$ then \hyperlink{om1}{$(\omega_1)$} holds true for $\omega_M$ if and only if
\begin{equation}\label{om1omegaMchar}
\exists\;L\in\NN_{>0}:\;\;\;\liminf_{j\rightarrow+\infty}\frac{(M_{Lj})^{1/(Lj)}}{(M_j)^{1/j}}>1.
\end{equation}
By definition \eqref{omegafromv} it is immediate that $\omega^v$ has \hyperlink{om1}{$(\omega_1)$} if and only if
\begin{equation}\label{om1forv}
\exists\;L\ge 1\;\forall\;t\ge 0:\;\;\;v^L(t)\le e^Lv_2(t).
\end{equation}

Summarizing, we get the following characterization (analogously to Lemma \ref{om6rem}):

\begin{lemma}\label{om1rem}
Let $u$ be a normalized and convex weight function. Then the following are equivalent:
\begin{itemize}
\item[$(a)$] $u$ satisfies \eqref{om1forv}.

\item[$(b)$] $\omega^u$ has \hyperlink{om1}{$(\omega_1)$}.

\item[$(c)$] $\omega_{M^u}$ has \hyperlink{om1}{$(\omega_1)$}.

\item[$(d)$] $M^u$ satisfies \eqref{om1omegaMchar}.
\end{itemize}
\end{lemma}

\subsection{Exponential-type weight systems in the literature}\label{refereecommentsect}
In \cite{Meise85}, \cite{Braun87}, \cite{MeiseTaylor87} the authors have already investigated weighted classes of entire functions defined in terms of (exponential-type) weight systems and, as examples, also weights defined in terms of a given sequence $M$ have been considered. More precisely, in these works the authors have focused on the study of (closed) ideals in such weighted spaces and also obtained sequence space representations; we comment now in detail:

\begin{itemize}
\item[$(a)$] In \cite{Meise85} the following setting is considered; see \cite[Def. 2.1-2.3]{Meise85}:
\begin{itemize}
\item[$(i)$] Weights $p$, not necessarily being radial, and weight systems $\mathbb{P}=(p_k)_k$ have been studied even on $\CC^n$. The notation $A_p(\CC^n)$ is used for $H^{\infty}_{\underline{\mathcal{V}}^{\mathfrak{c}}}$ when $\mathbb{P}:=(kp)_{k\in\NN_{>0}}$. Here $p$ is corresponding to $\omega^v$ (see \eqref{omegafromv}) and note: $(1)$ in Def. 2.1 is formally weaker than assuming that $\mathbb{P}$ consists of a non-decreasing sequence of weights whereas $(2)$ in Def. 2.1 is weaker than \hyperlink{om3}{$(\omega_3)$}. $(2)$ in Def. 2.2 for the system $(kp)_{k\in\NN_{>0}}$, when $p$ is radial, precisely means
    \begin{equation}\label{om6quasi}
    \forall\;k\in\NN_{>0}\;\exists\;m\in\NN_{>0}\;\exists\;L\ge 0\;\forall\;t\ge 0:\;\;\;2p(t)\le\frac{m}{k}p(t)+L,
    \end{equation}
which is clearly preserved under \hyperlink{sim}{$\sim$} and should be compared with \hyperlink{om6}{$(\omega_6)$} for $p$. For general systems $\mathbb{P}$ property \cite[Def. 2.2 $(2)$]{Meise85} can be viewed as a mixed variant of \eqref{om6quasi}.

But also \hyperlink{om1}{$(\omega_1)$} is appearing as an additional technical requirement; see e.g. \cite[Rem. 3.2, Thm. 4.7, Cor. 4.10, Prop. 4.15]{Meise85}.

\item[$(ii)$] In \cite[Sect. 2.6]{Meise85} weights $p=p_M=\omega_M$, and so expressed in terms of $M$, have been investigated; note that the basic assumption on $M$ there precisely corresponds to Definition \ref{defweightsequ}. Requirement \hyperlink{mg}{$(\on{mg})$} on $M$ is standard in order to ensure \hyperlink{om6}{$(\omega_6)$} for $\omega_M$; see \cite[Prop. 3.6]{Komatsu73}. Moreover, in \cite[Sect. 2.6 $(2)$]{Meise85} it has been mentioned that in a private communication H.-J. Petzsche has been able to show the statements \cite[Thm. 3.1, Prop. 3.4]{subaddlike} and this has also been mentioned in \cite[Rem. 8.9]{BraunMeiseTaylor90}. However, it seems that a proof has never been published by Petzsche and note that $\left(\frac{M_{jk}}{M_j^k}\right)^{1/j}=\left(\frac{M_{jk}^{1/k}}{M_j}\right)^{k/j}>1\Leftrightarrow\frac{(M_{jk})^{1/(jk)}}{(M_j)^{1/j}}>1$.

In \cite[Sect. 2.6 $(3)$, $(4)$]{Meise85} well-known examples have been mentioned: The first part $(3)$ deals with {\itshape Gevrey sequences} and the second one $(4)$ is corresponding to the sequences considered in \cite[Sect. 5.5]{whitneyextensionweightmatrix}. (For $s=2$ one gets the so-called \emph{$q$-Gevrey sequences} given by $(q^{j^2})_{j\in\NN}$, $q>1$.)

\item[$(iii)$] Next let us comment on the last sentence in \cite[Sect. 2.6 $(2)$]{Meise85} on p. 69 mentioning that $p_M$ is a radial weight on $\CC^n$: If $p:[0,+\infty)\rightarrow[0,+\infty)$ is non-decreasing and satisfies {\itshape both} \hyperlink{om6}{$(\omega_6)$} and \hyperlink{om1}{$(\omega_1)$}, i.e. \cite[Cor. 2.9 $(\ast)$]{Meise85}, then by iterating \hyperlink{om1}{$(\omega_1)$} one gets
$$\exists\;A,B\ge 1\;\forall\;t\ge 0:\;\;\;2p(t)\le Ap(t)+B,$$
hence \eqref{om6quasi} is valid.

We remark that the notation $M^k$ in \cite[Prop. 2.11]{Meise85} is not referring to the $k$-th power of $M$ but should be
viewed as a parameter/index; see the proof in \cite[Cor. 2.12]{Meise85}.

\item[$(iv)$] Finally, we comment on \cite[Sect. 5]{Meise85}: There a connection between non-increasing and non-decreasing exponential-type weight systems given in terms of special weight sequences is shown; see \cite[Prop. 5.3]{Meise85}. More precisely, the systems are expressed in terms of weight sequences $M$ and $N$ (resp. $M^k$, $N^k$, $k\in\NN_{>0}$) both (resp. for all $k$) having \hyperlink{mg}{$(\on{mg})$} and \eqref{om1omegaMchar} and such that $MN\hyperlink{approx}{\approx}(j!)_{j\in\NN}$ (resp. $M^kN^k\hyperlink{approx}{\approx}(j!)_{j\in\NN}$ for all $k$). This relation precisely means that, up to equivalence, the pair $M$ and $N$ (resp. for each $k$ the pair $M^k$ and $N^k$) is corresponding to the relation between $M$ and its so-called {\itshape conjugate sequence} $M^{*}$; see \cite[Sect. 2.5, Sect. 2.6]{microclasses}.
\end{itemize}

\item[$(b)$] In \cite{MeiseTaylor87} the analogous framework has been considered for the projective type:
\begin{itemize}
\item[$(i)$] $A^0_p$ denotes $H^{\infty}_{\overline{\mathcal{V}}^{\mathfrak{c}}}$ with the convention $p=\mathbb{P}:=(k^{-1}p)_{k\in\NN_{>0}}$. Here $p$ is corresponding to $\omega^v$ and \cite[Def. 1.1 $(2)$]{MeiseTaylor87} precisely denotes \hyperlink{om3}{$(\omega_3)$} and \cite[Def. 1.2 $(2)$]{MeiseTaylor87} is \eqref{om6quasi} with swapped order of quantifiers for $m$ and $k$. However, comment $(a)(iii)$ above applies also to this situation; see assumption \cite[Cor. 1.16 $(\ast)$]{MeiseTaylor87}.

\item[$(ii)$] And again \hyperlink{om1}{$(\omega_1)$} is appearing frequently as additional growth requirement. Note that \hyperlink{om1}{$(\omega_1)$} is precisely \cite[Thm. 2.7 $(1)$]{MeiseTaylor87} for $\mathbb{P}:=(k^{-1}p)_{k\in\NN_{>0}}$; see \cite[Cor. 2.8, Cor. 3.5]{MeiseTaylor87}.

\item[$(iii)$] In \cite[Ex. 3.7 $(3)$]{MeiseTaylor87} the weight sequence case corresponding to \cite[Sect. 2.6]{Meise85} has been revisited.

\item[$(iv)$] Finally, we mention that the assumptions made in \cite[Def. 4.1]{MeiseTaylor87} on the abstractly given system $Q=(q_k)_k$ seem to be not ''well-related'' to the weight sequence setting: Even convexity for each $q_k$ and $\lim_{t\rightarrow+\infty}\frac{q_k(t)}{t}=+\infty$ for all $k$ is required but both assumptions fail for $\omega_M$ and it seems that $q_k$ should be replaced by $k\varphi_{\omega_M}$ with $\varphi_{\omega_M}:=\omega_M\circ\exp$. However, then condition \cite[Def. 4.1 $(3)$]{MeiseTaylor87} differs from \hyperlink{om1}{$(\omega_1)$} for $\omega_M$ and also the definition of $A_q(\CC)$, corresponding to the system $(kq)_{k\in\NN}$, does not coincide with $H^{\infty}_{\underline{\mathcal{V}}^{\mathfrak{c}}}$. On the other hand, the results in this section seem to be analogous resp. dual to the ones shown in \cite[Sect. 5]{Meise85}; see comment $(a)(iv)$.
\end{itemize}

\item[$(c)$] In \cite{Braun87}, see the main results \cite[Thm. 2.7, Thm. 2.8]{Braun87}, the author has considered radial weights $p$ and used (again) the notation $A_p=H^{\infty}_{\underline{\mathcal{V}}^{\mathfrak{c}}}$, $A^0_p=H^{\infty}_{\overline{\mathcal{V}}^{\mathfrak{c}}}$. Also here $p$ is corresponding to $\omega^v$. Therefore, as it has been remarked by the referee, concerning the assumptions on $p$ in \cite[Sect. 1, p. 441]{Braun87} we note: subharmonicity in $(1)$ corresponds to \eqref{vconvexity}, $(3)$ is precisely \hyperlink{om3}{$(\omega_3)$} and $(4)$ is \hyperlink{om1}{$(\omega_1)$} for $\omega^v$. Thus the author deals with convex weights having \hyperlink{om1}{$(\omega_1)$}.
\end{itemize}

\section{Characterization of inclusions and equalities of weighted classes}\label{equsect}
In this section we give a complete characterization of the inclusion (and equality) of weighted classes of entire functions in terms of the growth of the corresponding weights. We treat both dilatation- and exponential-type weight systems. But before we focus on single (essential) weights and show that in this case one has an equivalence in \eqref{principalequ}.

\subsection{Characterization for essential weights}\label{essentialweightsect}
We recall the notions from \cite[Def. 1.1, Sect. 1.B]{BierstedtBonetTaskinen98}; see also \cite[Sect. 3]{Bonet2022survey}. Let a weight $v$ be given, then consider the function (radial growth condition) $w_v: \CC\rightarrow(0,+\infty)$ given by
\begin{equation}\label{corrgrowthcond}
w_v:=\frac{1}{v}
\end{equation}
and set
$$B_{w_v}(\CC):=\{f\in H(\CC): |f(z)|\le w_v(|z|),\;\forall\;z\in\CC\}.$$
The function $\widetilde{w_v}:\CC\rightarrow[0,+\infty)$ {\itshape associated with} $w_v$ is defined by
\begin{equation}\label{assoweightw}
\widetilde{w_v}(z):=\sup\{|f(z)|: f\in B_{w_v}(\CC)\},\;\;\;z\in\CC,
\end{equation}
and \cite[Observation 1.5]{BierstedtBonetTaskinen98} gives that $\widetilde{w_v}$ is again radial. Finally, set
\begin{equation}\label{assoweightv}
\widetilde{v}:=\frac{1}{\widetilde{w_v}}.
\end{equation}
If $v=v_M$, then we write $\widetilde{v_M}$. Similarly these notions can be introduced for arbitrary (open and connected) sets $G$ and even when using a more general notion for $v$ being a weight; e.g. see \cite{BierstedtBonetTaskinen98}, \cite{BonetDomanskiLindstroemTaskinen98} ($v$ is only assumed to be continuous and strictly positive). In \cite[Prop. 1.3]{BonetDomanskiLindstroemTaskinen98} it has been shown that isometrically $H^{\infty}_{v}(\mathbb{D})=H^{\infty}_{\widetilde{v}}(\mathbb{D})$ but the proof there gives $H^{\infty}_{v}(\CC)=H^{\infty}_{\widetilde{v}}(\CC)$ too; see also \cite[Observ. 1.12]{BierstedtBonetTaskinen98}.

Finally, we have the following definition; see \cite[Sect. 3]{Bonet2022survey}.

\begin{definition}
A weight $v$ is called {\itshape essential} if $v\hyperlink{sim}{\sim}\widetilde{v}$.
\end{definition}

We can now formulate the main theorem in this context. It directly follows from the more general results \cite[Prop. 2.1, Cor. 2.2]{BonetDomanskiLindstroemTaskinen98} dealing with the composition operator $C_{\varphi}(f):=f\circ\varphi$ acting on weighted classes of analytic functions defined on $\mathbb{D}$.

\begin{theorem}\label{essentialweightcharactnew}
Let $u$ and $v$ be weights. Then the following are equivalent:
\begin{itemize}
\item[$(i)$] The weights are related by $\widetilde{u}\hyperlink{ompreceq}{\preceq}v$.

\item[$(ii)$] The weights are related by $\widetilde{u}\hyperlink{ompreceq}{\preceq}\widetilde{v}$.

\item[$(iii)$] We have $H^{\infty}_u(G)\subseteq H^{\infty}_v(G)$ (with continuous inclusion) for $G=\mathbb{D}$ and/or $G=\CC$.
\end{itemize}
Thus, if $u$ is essential, then the inclusion in $(iii)$ holds if and only if $u\hyperlink{ompreceq}{\preceq}v$.
\end{theorem}

When $G=\mathbb{D}$ then the relation $u\hyperlink{ompreceq}{\preceq}v$ means $v(t)\le C u(t)$ for all $0\le t<1$ and some $C\ge 1$.

\demo{Proof}
We apply \cite[Prop. 2.1]{BonetDomanskiLindstroemTaskinen98} to $\varphi=\id$ and hence $C_{\id}=\id$; see also \cite[Cor. 2.2]{BonetDomanskiLindstroemTaskinen98}. Note that in \cite[Prop. 2.1]{BonetDomanskiLindstroemTaskinen98} the authors deal with $G=\mathbb{D}$ but the proof there is also valid for $G=\CC$.
\qed\enddemo

However, we give a second independent proof of the previous result by involving different techniques (optimal functions).

\begin{theorem}\label{essentialweightcharact}
Let $u$ and $v$ be weights and assume that $u$ is essential. Then the following are equivalent:
\begin{itemize}
\item[$(i)$] The weights are related by $u\hyperlink{ompreceq}{\preceq}v$.

\item[$(ii)$] We have $H^{\infty}_u(G)\subseteq H^{\infty}_v(G)$ (with continuous inclusion) for $G=\mathbb{D}$ and/or $G=\CC$.
\end{itemize}
\end{theorem}

\demo{Proof}
It remains to prove $(ii)\Rightarrow(i)$. For $G=\CC$ we use \cite[Thm. 2]{AbakumovDoubtsov18} (see also \cite[Thm. 3]{Bonet2022survey}) and get that $u$ is essential if and only if $u$ is {\itshape approximable by the maximum of an analytic function modulus,} i.e.
\begin{equation}\label{approximable}
\exists\;f_u\in H(\CC):\;\;\;\frac{1}{u}\hyperlink{sim}{\sim} t\mapsto M(f_u,t),
\end{equation}
with $M(f,t):=\sup\{|f(z)|: |z|=t\}$; see \cite[Def. 3]{AbakumovDoubtsov18} and \cite[p. 3, p. 5]{Bonet2022survey}. Recall that $w=\frac{1}{u}$ denotes the weight used in the notation in \cite{AbakumovDoubtsov18}. \eqref{approximable} gives that $f_u\in H^{\infty}_u(\CC)$ and more precisely this relation ensures the existence of {\itshape optimal functions} $f_u\in H^{\infty}_u(\CC)$ admitting sharp estimates from below.

Then \eqref{approximable} and the inclusion $H^{\infty}_u(\CC)\subseteq H^{\infty}_v(\CC)$ (as sets) give
$$\exists\;C,D\ge 1\;\forall\;t\ge 0:\;\;\;\frac{1}{u(t)}\le CM(f_u,t)\le CD\frac{1}{v(t)},$$
i.e. $v(t)=O(u(t))$ and so $u\hyperlink{ompreceq}{\preceq}v$ is verified.\vspace{6pt}

When $G=\mathbb{D}$, then by \cite[Thm. 1.2]{AbakumovDoubtsov15}, see also \cite[Thm. 1, Rem. 1, Def. 1, Thm. 5]{AbakumovDoubtsov18} and \cite[Thm. 2]{Bonet2022survey}, again $u$ is essential if and only if \eqref{approximable} is valid. The rest follows as before.

Note that for this case in order to apply the cited results we do not require necessarily that $u$ is rapidly decreasing.
\qed\enddemo

The next result illustrates that the notion of essentiality and the weight sequence setting are very closely connected via involving the associated weight sequence.

\begin{theorem}\label{essentialweightsequthm}
We get the following:
\begin{itemize}
\item[$(i)$] Let $M$ be a log-convex weight sequence. Then $v_M$ is essential (for the weighted entire case) and in fact we even get $v_M=\widetilde{v_M}$.

\item[$(ii)$] The following are equivalent:
\begin{itemize}
\item[$(a)$] $u$ is an essential weight.

\item[$(b)$] $u\hyperlink{sim}{\sim}v_{M^u}$ is valid.
\end{itemize}
\end{itemize}
\end{theorem}

By recalling \eqref{powermultisequequ} and \eqref{powermultisequequ1} we thus have that $v_{M,c}$ is essential and $\widetilde{v_{M,c}}=v_{M,c}$ for each $c>0$. The difference between $u\hyperlink{sim}{\sim}v_{M^u}$ and the weaker relation \eqref{omegavequiv} (due to the appearing quadratic power) illustrates that (log)-convexity is in general a weaker requirement than essentiality; see also Remark \ref{AtoEremark} for more explanations.

\demo{Proof}
$(i)$ For any weight $u$ in \cite[Lemma 1]{AbakumovDoubtsov18}, by using a statement from \cite{ErdoesKoevari1956}, it has been shown that (see also \cite[p.5]{Bonet2022survey}):
\begin{equation}\label{ErdoesKoevconsequence}
\forall\;t\ge 0:\;\;\;P_u(t)\le\frac{1}{\widetilde{u}(t)}\le 6P_u(t).
\end{equation}
For this recall the notation $w=\frac{1}{u}$, $\widetilde{w}=\frac{1}{\widetilde{u}}$ used in \cite{AbakumovDoubtsov18}. Therefore, by combining \eqref{ErdoesKoevconsequence} with Lemma \ref{monomialnormlemma} applied to $u=v_M$ we get
$$\forall\;t\ge 0:\;\;\;\frac{1}{v_M(t)}\le\frac{1}{\widetilde{v_M}(t)}\le\frac{6}{v_M(t)},$$
i.e. $v_M\hyperlink{sim}{\sim}\widetilde{v_M}$ is verified and so $v_M$ is essential.

Moreover, for any growth condition $w$ we get $0\le\widetilde{w}\le w$; see \cite[Prop. 1.2]{BierstedtBonetTaskinen98}. Thus by definitions \eqref{corrgrowthcond}, \eqref{assoweightv} and the first part in \eqref{ErdoesKoevconsequence} $$v_M=\frac{1}{w_{v_M}}\le\frac{1}{\widetilde{w_{v_M}}}=\widetilde{v_M}\le v_M$$
holds and we are done.\vspace{6pt}

$(ii)(a)\Rightarrow(b)$ Since $u$ is essential the non-trivial estimate $\widetilde{u}(t)\le Cu(t)$ holds true for some $C\ge 1$ and all $t\ge 0$. Then the second part in \eqref{ErdoesKoevconsequence} and \eqref{monomialnormlemmaequ} yield
$$\forall\;t\ge 0:\;\;\;\frac{1}{v_{M^u}(t)}=P_u(t)\le\frac{1}{u(t)}\le\frac{C}{\widetilde{u}(t)}\le 6CP_u(t)=\frac{6C}{v_{M^u}(t)},$$
and so $v_{M^u}\hyperlink{sim}{\sim}u$ is verified.

$(ii)(b)\Rightarrow(a)$ The assumption and the first part in \eqref{ErdoesKoevconsequence} give
$$\exists\;C\ge 1\;\forall\;t\ge 0:\;\;\;\frac{1}{Cu(t)}\le\frac{1}{v_{M^u}(t)}=P_u(t)\le\frac{1}{\widetilde{u}(t)},$$
hence $\widetilde{u}(t)\le Cu(t)$. Since always $u(t)\le\widetilde{u}(t)$ holds true (see \cite[Prop. 1.2]{BierstedtBonetTaskinen98}) the relation $u\hyperlink{sim}{\sim}\widetilde{u}$ is verified and we are done.
\qed\enddemo

Finally, when combining Theorems \ref{essentialweightsequthm} and \ref{essentialweightcharact} with Proposition \ref{strongNMrelationlemma} we arrive at the following characterization:

\begin{theorem}\label{weightholombysequcharactsingle}
Let $M$ and $N$ be log-convex weight sequences. Then both $v_M$ and $v_N$ are essential and the following are equivalent:
\begin{itemize}
\item[$(i)$] The sequences satisfy \eqref{strongNMrelation}, i.e.
$$\exists\;A\ge 1\;\forall\;j\in\NN:\;\;\;N_j\le AM_j.$$

\item[$(ii)$] The weights satisfy
$$v_M\hyperlink{ompreceq}{\preceq}v_N.$$

\item[$(iii)$] We have
$$H^{\infty}_{v_M}(\CC)\subseteq H^{\infty}_{v_N}(\CC),$$
with continuous inclusion.
\end{itemize}
Thus \eqref{LCapproxremequ} if and only if $v_M\hyperlink{sim}{\sim}v_N$ if and only if $H^{\infty}_{v_M}(\CC)=H^{\infty}_{v_N}(\CC)$ (as Banach spaces).
\end{theorem}

\begin{remark}\label{AtoEremark}
\emph{We comment on the known relation between essentiality and convexity and its consequences.}

\begin{itemize}
\item[$(A)$] \emph{Let $u$ and $r$ be given weights. Of course, Theorem \ref{weightholombysequcharactsingle} applies to the associated sequences $M^u$ and $M^r$ and the proof of Theorem \ref{essentialweightsequthm} yields $v_{M^u}\hyperlink{sim}{\sim}\widetilde{u}$ and $v_{M^r}\hyperlink{sim}{\sim}\widetilde{r}$. The first part in \eqref{omegavequivnewnew}, for which even convexity is not required, gives the following: Assume that the inclusion $H^{\infty}_{u}(\CC)\subseteq H^{\infty}_{r}(\CC)$ is valid (as sets), then $H^{\infty}_{v_{M^u}}(\CC)\subseteq H^{\infty}_{r}(\CC)$ and by Theorem \ref{essentialweightsequthm} the weight $v_{M^u}$ is essential. Hence $v_{M^u}\hyperlink{ompreceq}{\preceq}r$ follows by Theorem \ref{essentialweightcharact}.}

\emph{However, from this the desired relation $u\hyperlink{ompreceq}{\preceq}r$ does not follow automatically which is due to the quadratic power appearing in \eqref{omegavequiv}. In view of $(ii)$ in Theorem \ref{essentialweightsequthm} we can conclude provided that $u$ is essential.}

\item[$(B)$] \emph{By \cite{AbakumovDoubtsov15}, \cite{AbakumovDoubtsov18} for the weighted entire case each essential weight $u$ is equivalent to a convex one. Since $v_{M^u}$ is convex, $(ii)$ in Theorem \ref{essentialweightsequthm} proves this implication and provides more precise information.}

   \emph{ For $G=\mathbb{D}$ one even has that $u$ is essential, equivalently \eqref{approximable} holds true, if and only if it is equivalent to a convex weight.}

\item[$(C)$] \emph{In the discussion after \cite[Thm. 1]{AbakumovDoubtsov18} resp. after \cite[Thm. 2]{Bonet2022survey} it is explained that if $u$ is rapidly decreasing, in particular for weight functions in our setting, in general the characterization \cite[Thm. 2]{Bonet2022survey} fails for $G=\CC$. More precisely, there exist (convex!) weights $u$ on $\CC$ such that \eqref{approximable} is violated. In particular, by \cite[Thm. 2]{AbakumovDoubtsov18}, \cite[Thm. 3]{Bonet2022survey}, such weights cannot be essential and by $(ii)$ in Theorem \ref{essentialweightsequthm} they cannot be equivalent to $v_{M^u}$. In any case, for non-essential weights $u$ defined on $\CC$ property \eqref{approximable} fails.}

\item[$(D)$] \emph{Hence by the previous comments we see that for the weighted entire case when dealing with general convex weights the proof of $\Leftarrow$ in \eqref{principalequ} does not follow from (the proof of) Theorem \ref{essentialweightcharact} since \eqref{approximable} fails in general. On the other hand convexity of $u$ is sufficient for our techniques (for involving the sequence $M^u$; see Theorem \ref{weightholombyfctcharact}). Note that in view of Remark \ref{assoweightsequ} the assumption that $u$ is rapidly decreasing should be considered as standard when involving $M^u$.}

\item[$(E)$] \emph{In this context also note that \eqref{approximable} implies \eqref{Clunieremarkequ} (with $\phi=-\log\circ u$ there); see Remark \ref{Clunieremark} for more explanations.}

\item[$(F$)] \emph{In \cite[Prop. 3.1]{BierstedtBonetTaskinen98} sufficient conditions on (log)-convex weights $v$ are given to ensure essentiality.}
\end{itemize}
\end{remark}

We close this section by showing a partial converse to $(B)$.

\begin{theorem}\label{assoweightfctprop}
Let $u$ be a normalized and convex weight function.

\begin{itemize}
\item[$(i)$] We get
\begin{equation}\label{assoweightfctpropequ0}
\exists\;B\ge 1\;\forall\;t\ge 0:\;\;\;u(t)\le\widetilde{u}(t)\le Bu^{1/2}(t),
\end{equation}
i.e. $\widetilde{u}\hyperlink{ompreceq}{\preceq}u\hyperlink{ompreceqpowc}{\preceq^{\mathfrak{c}}}\widetilde{u}$. Thus $u\hyperlink{simpowc}{\sim^{\mathfrak{c}}}\widetilde{u}$ holds true.

\item[$(ii)$] If in addition $u$ is of moderate growth, then
\begin{equation}\label{assoweightfctpropequ1}
\exists\;B,H\ge 1\;\forall\;t\ge 0:\;\;\;u(t)\le\widetilde{u}(t)\le Bu_{H^{-1}}(t)=Bu(t/H),
\end{equation}
i.e. $\widetilde{u}\hyperlink{ompreceq}{\preceq}u\hyperlink{ompreceqc}{\preceq_{\mathfrak{c}}}\widetilde{u}$. Thus $u\hyperlink{simc}{\sim_{\mathfrak{c}}}\widetilde{u}$ holds true.
\end{itemize}
\end{theorem}

\demo{Proof}
Let us write $w_u$ and $w_{M,c}:=w_{v_{M,c}}$ for the corresponding growth conditions; see \eqref{corrgrowthcond}. We will also use the following immediate properties; see e.g. \cite[Prop. 1.2]{BierstedtBonetTaskinen98}: For any growth condition $w$ we get $0\le\widetilde{w}\le w$, write $(\star)$ for this, and given weights $u,v$ with $u\le v\Leftrightarrow w_v\le w_u$, then $\widetilde{w_v}\le\widetilde{w_u}$, written as relation $(\star\star)$.\vspace{6pt}

$(i)$ First, \eqref{omegavequivnew} applied to the weight $u$ transfers into
\begin{equation}\label{omegavequivforw}
\exists\;A\ge 1\;\forall\;c>0\;\forall\;t\ge 0:\;\;\;w_{M^u,c}(t)\le w_{u_c}(t)\le A w_{M^u,c}(t)^2.
\end{equation}
Take $c=1$ and from this we get
\begin{equation}\label{omegavequivforw0}
\exists\;A\ge 1\;\forall\;t\ge 0:\;\;\;\frac{1}{\sqrt{A}}\sqrt{w_u(t)}\le w_{M^u}(t)=\widetilde{w_{M^u}}(t)\le\widetilde{w_u}(t).
\end{equation}
More precisely, the first estimate holds by the second half in \eqref{omegavequivforw}, the second estimate since $v_{M^u}$ is essential and coincides with $\widetilde{v_{M^u}}$ by Theorem \ref{essentialweightsequthm}. The third one holds by the first estimate in \eqref{omegavequivforw} and $(\star\star)$.

So \eqref{omegavequivforw0} is verified and thus the second half of \eqref{assoweightfctpropequ0} follows by \eqref{corrgrowthcond}, \eqref{assoweightv} (the first half is precisely $(\star)$).\vspace{6pt}

$(ii)$ In view of \eqref{corrgrowthcond} and since for any weight $u$ and $c>0$ we get
\begin{equation}\label{assoweightfctpropaux}
w_{u_c}(t)=\frac{1}{u_c(t)}=\frac{1}{u(ct)}=w_u(ct)=(w_u)_{c}(t),
\end{equation}
a weight $u$ is of moderate growth if and only if the corresponding growth condition satisfies (with the same $H$)
\begin{equation}\label{om6forw}
\exists\;H\ge 1\;\forall\;t\ge 0:\;\;\;w_u(t)^2\le e^Hw_{u_H}(t)=e^Hw_u(Ht).
\end{equation}
Next let us verify
\begin{equation}\label{assoweightfctpropaux0}
\exists\;A\;\ge 1\;\forall\;c>0\;\forall\;t\ge 0:\;\;\;w_{u_c}(t)\le A w_{M^u,c}(t)^2=A\widetilde{w_{M^u,c}}(t)^2\le A\widetilde{w_{u_c}}(t)^2.
\end{equation}
The first estimate holds by \eqref{omegavequivforw} (second part), the second one since each $v_{M^u,c}$ is essential by Theorem \ref{essentialweightsequthm} and finally the third estimate in \eqref{assoweightfctpropaux0} holds by combining the first part of \eqref{omegavequivforw} with $(\star\star)$.

Now we continue \eqref{assoweightfctpropaux0} (with $c=H^{-1}$) and prove
\begin{equation}\label{assoweightfctpropaux1}
A\widetilde{w_{u_{H^{-1}}}}(t)^2\le Ae^H\widetilde{w_u}(t).
\end{equation}
By \eqref{om6forw} one has $e^Hw_{u}(t)\ge w_{u_{H^{-1}}}(t)^2$ and then note that for any $f\in B_w(\CC)$ and $c\in\NN_{>0}$ we have $f^c\in B_{w^c}(\CC)$. So
$$\widetilde{w}^c(z)=\sup\{|f(z)|^c: f\in B_w(\CC)\}=\sup\{|f^c(z)|: f\in B_w(\CC)\}\le\sup\{|g(z)|: g\in B_{w^c}(\CC)\}=\widetilde{w^c}(z)$$
follows and by combining both estimates (take $c=2$) with $(\star\star)$ and the fact that $\widetilde{Cw}=C\widetilde{w}$ for any $C>0$ and any growth condition $w$ (see again \cite[Prop. 1.2]{BierstedtBonetTaskinen98}) the estimate \eqref{assoweightfctpropaux1} holds.\vspace{6pt}

Summarizing, by combining \eqref{assoweightfctpropaux0} and \eqref{assoweightfctpropaux1} we have verified
$$\exists\;A,H\ge 1\;\forall\;t\ge 0:\;\;\;w_{u_{H^{-1}}}(t)\le Ae^H\widetilde{w_u}(t),$$
and so by taking into account \eqref{corrgrowthcond}, \eqref{assoweightv} we get the second half of \eqref{assoweightfctpropequ1} with $B:=Ae^H$ and $H$ being the constant appearing in \eqref{om6forv} for $u$.
\qed\enddemo

\subsection{The dilatation-type weight sequence setting}\label{equsectweightsequdila}
We start with the following:

\begin{lemma}\label{charactprop0}
Let $M$ and $N$ be weight sequences and assume that $N$ is log-convex. Consider the following assertions:
\begin{itemize}
	\item[$(i)$] We have \eqref{strongNMrelation}; i.e.
$$\exists\;A\ge 1\;\forall\;j\in\NN:\;\;\;N_j\le AM_j.$$

\item[$(ii)$] The weights satisfy
\begin{equation}\label{forallcpreceq}
	\forall\;c>0:\;\;\;v_{M,c}\hyperlink{ompreceq}{\preceq}v_{N,c},
\end{equation}
which implies $v_M\hyperlink{ompreceqc}{\preceq_{\mathfrak{c}}}v_N$.

\item[$(iii)$] We have that
\begin{equation}\label{forallcpreceq1}
\forall\;c>0:\;\;\;H^{\infty}_{v_{M,c}}\subseteq H^{\infty}_{v_{N,c}},
\end{equation}
with continuous inclusion.
\end{itemize}
Then $(i)\Leftrightarrow(ii)\Rightarrow(iii)$ holds true. Alternatively, in \eqref{forallcpreceq} and \eqref{forallcpreceq1} we can assume $\exists\;c>0$ instead and $(iii)$ implies both $H^{\infty}_{\underline{\mathcal{M}}_{\mathfrak{c}}}\subseteq H^{\infty}_{\underline{\mathcal{N}}_{\mathfrak{c}}}$ and $H^{\infty}_{\overline{\mathcal{M}}_{\mathfrak{c}}}\subseteq H^{\infty}_{\overline{\mathcal{N}}_{\mathfrak{c}}}$ with continuous inclusions.
\end{lemma}

If log-convexity for $N$ fails, then only $(i)\Rightarrow(ii)\Rightarrow(iii)$ is valid.

\demo{Proof}
$(i)\Leftrightarrow(ii)$ follows from resp. as in Proposition \ref{strongNMrelationlemma} and $(ii)\Rightarrow(iii)$ is clear by definition \eqref{weightedclassdef}.
\qed\enddemo

In the next result we replace \eqref{strongNMrelation} by the weaker requirement $N\hyperlink{preceq}{\preceq}M$.

\begin{proposition}\label{charactprop}
Let $M$ and $N$ be weight sequences. Consider the following assertions:
\begin{itemize}
\item[$(i)$] $N\hyperlink{preceq}{\preceq}M$,
\item[$(ii)$] $v_M\hyperlink{ompreceqc}{\preceq_{\mathfrak{c}}}v_N$,
\item[$(iii)$]
$$H^{\infty}_{\underline{\mathcal{M}}_{\mathfrak{c}}}\subseteq H^{\infty}_{\underline{\mathcal{N}}_{\mathfrak{c}}},\hspace{15pt}\text{and}\hspace{15pt}H^{\infty}_{\overline{\mathcal{M}}_{\mathfrak{c}}}\subseteq H^{\infty}_{\overline{\mathcal{N}}_{\mathfrak{c}}},$$
with continuous inclusions.
\end{itemize}
Then $(i)\Rightarrow(ii)\Rightarrow(iii)$ is valid.
\end{proposition}

\demo{Proof}
$(i)\Rightarrow(ii)$ By $N\hyperlink{preceq}{\preceq}M$ there exist $C,h\ge 1$ such that $N_j\le Ch^jM_j$ for all $j\in\NN$ and by definition of associated weight functions we get
\begin{equation}\label{assoweightsrelated}
\exists\;C,h\ge 1\;\forall\;t\ge 0:\;\;\;\omega_M(t)\le\omega_N(ht)+\log(C),
\end{equation}
hence in view of \eqref{weights}
\begin{equation}\label{assoweightsrelatedforv}
\exists\;h\ge 1\;\forall\;c>0:\;\;\;v_{N,hc}(t)=O(v_{M,c}(t)),\;\;\;t\rightarrow+\infty.
\end{equation}
This is precisely $v_M\hyperlink{ompreceqc}{\preceq_{\mathfrak{c}}}v_N$.\vspace{6pt}

$(ii)\Rightarrow(iii)$ This is clear by definition (recall Proposition \ref{firstweighfctprop}).
\qed\enddemo

Now let us see to prove a complete characterization. In order to do so, in the inductive structure, the idea is to use an optimal function. For the weight sequence approach by the previous section we can use $f_{v_{M,c}}$ from \eqref{approximable}. We present here a different, weaker but sufficient construction involving a dilatation parameter. Let $M$ be a weight sequence, then consider the entire function
\begin{equation}\label{charholomfct}
\theta_{M,c}(z):=\sum_{j\ge 0}\frac{1}{2^jM_j}(cz)^j,\;\;\;z\in\CC,\;c>0.
\end{equation}
This definition is inspired by \cite{integralprescribed}, more precisely compare $\theta_{M,c}$ with $F$ in \cite{integralprescribed}; see also Remark \ref{Clunieremark}.

\begin{lemma}\label{charholomfctlemma}
Let $M$ be a weight sequence. Then
$$\forall\;c>0:\;\;\;\theta_{M,c}\in H^{\infty}_{v_{M,c}}\subseteq H^{\infty}_{\underline{\mathcal{M}}_{\mathfrak{c}}},$$
and
\begin{equation}\label{charholomfctlemmaequ}
\forall\;c>0\;\forall\;t\ge 0:\;\;\;\exp(\omega_M(ct/2))\le|\theta_{M,c}(t)|=\theta_{M,c}(t).
\end{equation}
\end{lemma}

\demo{Proof}
For any $z\in\CC$ we estimate by
$$|\theta_{M,c}(z)|\le\sum_{j\ge 0}\frac{1}{2^jM_j}(c|z|)^j\le\sup_{k\in\NN}\frac{(c|z|)^k}{M_k}\sum_{j\ge 0}\frac{1}{2^j}=2\exp(\omega_M(c|z|)),$$
showing $\theta_{M,c}\in H^{\infty}_{v_{M,c}}$. If $t\ge 0$, then $\theta_{M,c}(t)\in\RR_{\ge 0}$ and we get
$$\theta_{M,c}(t)=\sum_{j\ge 0}\frac{1}{2^jM_j}(ct)^j\ge\sup_{k\in\NN}\frac{(ct)^k}{2^kM_k}=\exp(\omega_M(ct/2)),$$
since each summand is non-negative.
\qed\enddemo

\begin{remark}\label{characteristicnotcontained}
\emph{The difference between the estimate \eqref{charholomfctlemmaequ} for $\theta_{M,c}$ and $f_{v_{M,c}}$ in \eqref{approximable} is that $|\theta_{M,c}|\hyperlink{sim}{\sim}\frac{1}{v_{M,c}}$ is not clear in general since an additional dilatation parameter $\frac{1}{2}$ appears on the left-hand side in \eqref{charholomfctlemmaequ}. However, it turns out that this fact is not disturbing for dilatation-type systems. On the other hand in the previous Lemma formally log-convexity for $M$ is not required but in view of Remark \ref{LCwlogrem} this difference is negligible.}\vspace{6pt}

\emph{Let $M$ be a log-convex weight sequence. Any entire function with the property that the estimate from \eqref{charholomfctlemmaequ} holds for some arbitrary but fixed parameter $c>0$ cannot be contained in $H^{\infty}_{\overline{\mathcal{M}}_{\mathfrak{c}}}$: Assume that there exists $\psi\in H^{\infty}_{\overline{\mathcal{M}}_{\mathfrak{c}}}$ such that \eqref{charholomfctlemmaequ} for some $c>0$, then
$$\forall\;d>0\;\exists\;D\ge 1\;\forall\;t\ge 0:\;\;\;\exp(\omega_M(ct/2))\le|\psi(t)|\le D\exp(\omega_M(dt)),$$
which implies $\omega_M(ct/2)-\omega_M(dt)\le\log(D)$. But this is contradicting \eqref{strongweightsrelation} as $t\rightarrow+\infty$ when $d<c/2$.}

\emph{In particular this comment also applies to $f_{v_{M,c}}$ from \eqref{approximable} and hence for projective structures a different technique is required.}
\end{remark}

\begin{remark}\label{poweroptimalrem}
\emph{Similarly, for any $c\in\NN_{>0}$ we can put}
\begin{equation}\label{charholomfctpow}
\theta^c_M(z):=\sum_{j\ge 0}\frac{1}{2^j(M_j)^c}z^{cj},\;\;\;z\in\CC,
\end{equation}
\emph{and show analogously as before that $\theta^c_M\in H^{\infty}_{v^c_M}\subseteq H^{\infty}_{\underline{\mathcal{M}}^{\mathfrak{c}}}$ and}
\begin{equation}\label{charholomfctlemmaequ1}
\forall\;c\in\NN_{>0}\;\forall\;t\ge 0:\;\;\;\exp(c\omega_M(t/2^{1/c}))\le|\theta^c_M(t)|=\theta^c_M(t).
\end{equation}
\emph{However, it turns out that for exponential-type structures this construction is less optimal than applying the functions from \eqref{approximable} since in \eqref{charholomfctlemmaequ1} the appearing dilatation parameter $\frac{1}{2^{1/c}}$ requires additional technical assumptions on the weights; see Section \ref{equsectweightsequexpo} and in particular Remark \ref{poweroptimalrem1} for more details.}
\end{remark}

\begin{proposition}\label{charactprop1}
Let $M,N$ be weight sequences. Assume that $N$ is log-convex and that (as sets)
\begin{equation}\label{charactprop1equ}
\exists\;c>0\;\exists\;d>0:\;\;\;H^{\infty}_{v_{M,c}}\subseteq H^{\infty}_{v_{N,d}}.
\end{equation}
Then $N\hyperlink{preceq}{\preceq}M$ follows.

Alternatively, in order to conclude we can replace in \eqref{charactprop1equ} one or both $\exists$ by $\forall$.
\end{proposition}

For the proof only log-convexity for $N$ is required but indispensable (by using \eqref{Prop32Komatsu}); otherwise we only get $N^{\on{lc}}\hyperlink{preceq}{\preceq}M^{\on{lc}}(\le M)$, i.e. the desired relation for the log-convex minorants.

\demo{Proof}
Let the parameters $c$ and $d$ be given. The inclusion $H^{\infty}_{v_{M,c}}\subseteq H^{\infty}_{v_{N,d}}$ and Lemma \ref{charholomfctlemma} yield $\theta_{M,c}\in H^{\infty}_{v_{N,d}}$ and so
$$\exists\;C\ge 1\;\forall\;t\ge 0:\;\;\;\exp(\omega_M(ct/2))\le\theta_{M,c}(t)=|\theta_{M,c}(t)|\le C\exp(\omega_N(dt)).$$
Hence, by applying \eqref{Prop32Komatsu} we obtain for all $j\in\NN$ that
\begin{align*}
M_j&\ge M^{\on{lc}}_j=\sup_{t\ge 0}\frac{t^j}{\exp(\omega_{M}(t))}=\sup_{s\ge 0}\frac{(cs/2)^j}{\exp(\omega_{M}(cs/2))}\ge\frac{(c/2)^j}{C}\sup_{s\ge 0}\frac{s^j}{\exp(\omega_{N}(ds))}
\\&
=\frac{(c/(2d))^j}{C}\sup_{u\ge 0}\frac{u^j}{\exp(\omega_{N}(u))}=\frac{(c/(2d))^j}{C}N_j,
\end{align*}
i.e. $N\hyperlink{preceq}{\preceq}M$ is verified.

A different technique is to apply Theorem \ref{weightholombysequcharactsingle} to $v_{M^{\frac{1}{c}}}$ and $v_{N^{\frac{1}{d}}}$ (see \eqref{powermultisequequ1}) and get
$\frac{1}{d^j}N_j=N^{\frac{1}{d}}_j\le AM^{\frac{1}{c}}_j=A\frac{1}{c^j}M_j$. Note that in order to apply Theorem \ref{weightholombysequcharactsingle} we have to assume that {\itshape both} sequences are log-convex.
\qed\enddemo

For the characterization of $H^{\infty}_{\overline{\mathcal{M}}_{\mathfrak{c}}}\subseteq H^{\infty}_{\overline{\mathcal{N}}_{\mathfrak{c}}}$ by Remark \ref{characteristicnotcontained} we have to come up with a different technique; we use the family of monomials $(f_j)_{j\in\NN}$.

\begin{proposition}\label{charactprop2}
Let $M,N$ be weight sequences. Assume that $N$ is log-convex and (with continuous inclusion)
\begin{equation}\label{charactprop2equ}
H^{\infty}_{\overline{\mathcal{M}}_{\mathfrak{c}}}\subseteq H^{\infty}_{\overline{\mathcal{N}}_{\mathfrak{c}}}.
\end{equation}
Then $N\hyperlink{preceq}{\preceq}M$ holds true.
\end{proposition}

\demo{Proof}
Both $H^{\infty}_{\overline{\mathcal{M}}_{\mathfrak{c}}}$ and $H^{\infty}_{\overline{\mathcal{N}}_{\mathfrak{c}}}$ are Fr\'{e}chet and by the closed graph theorem the inclusion in \eqref{charactprop2equ} is continuous since convergence in $H^{\infty}_{\overline{\mathcal{M}}_{\mathfrak{c}}}$ implies point-wise convergence; see the proof of \cite[Prop. 4.6 $(1)$]{compositionpaper} and more generally \cite[Prop. 4.5, Rem. 4.6]{PTTvsmatrix}.

So
$$\forall\;d>0\;\exists\;c>0\;\exists\;C\ge 1\;\forall\;f\in H^{\infty}_{\overline{\mathcal{M}}_{\mathfrak{c}}}:\;\;\; \|f\|_{v_{N,d}}\le C\|f\|_{v_{M,c}}.$$
We apply this estimate to the family of monomials $(f_k)_{k\in\NN}$ and get
$$\forall\;d>0\;\exists\;c>0\;\exists\;C\ge 1\;\forall\;k\in\NN\;\forall\;z\in\CC:\;\;\;\frac{|z|^k}{\exp(\omega_N(d|z|))}\le C\frac{|z|^k}{\exp(\omega_M(c|z|))}.$$
Thus, by using again \eqref{Prop32Komatsu}, we obtain for all $j\in\NN$ that
\begin{align*}
N_j&=\sup_{t\ge 0}\frac{t^j}{\exp(\omega_{N}(t))}=\sup_{s\ge 0}\frac{(ds)^j}{\exp(\omega_{N}(ds))}\le Cd^j\sup_{s\ge 0}\frac{s^j}{\exp(\omega_{M}(cs))}
\\&
=C\left(\frac{d}{c}\right)^j\sup_{u\ge 0}\frac{u^j}{\exp(\omega_{M}(u))}=C\left(\frac{d}{c}\right)^jM^{\on{lc}}_j\le C\left(\frac{d}{c}\right)^jM_j,
\end{align*}
i.e. $N\hyperlink{preceq}{\preceq}M$ is verified.
\qed\enddemo

Combining now Propositions \ref{charactprop}, \ref{charactprop1} and \ref{charactprop2} we arrive at the following characterization.

\begin{theorem}\label{weightholombysequcharact}
Let $M,N$ be weight sequences and assume that $N$ is log-convex. Then the following are equivalent:
\begin{itemize}
\item[$(a)$] The sequences satisfy $N\hyperlink{preceq}{\preceq}M$.

\item[$(b)$] The weights $v_M$ and $v_N$ satisfy $v_M\hyperlink{ompreceqc}{\preceq_{\mathfrak{c}}}v_N$.

\item[$(c)$] We have
$$H^{\infty}_{\underline{\mathcal{M}}_{\mathfrak{c}}}\subseteq H^{\infty}_{\underline{\mathcal{N}}_{\mathfrak{c}}},$$
with continuous inclusion.

\item[$(d)$] We have
$$H^{\infty}_{\overline{\mathcal{M}}_{\mathfrak{c}}}\subseteq H^{\infty}_{\overline{\mathcal{N}}_{\mathfrak{c}}},$$
    with continuous inclusion.
\end{itemize}
Consequently, relation \hyperlink{simc}{$\sim_{\mathfrak{c}}$} is characterizing the equivalence of the classes (as l.c.v.s.).
\end{theorem}

{\itshape Note:} When assuming instead of \eqref{charactprop2equ} the (continuous) inclusion $H^{\infty}_{v_{M,c}}\subseteq H^{\infty}_{v_{N,c}}$, $c>0$ arbitrary, and following the proof in Proposition \ref{charactprop2} with $c=d$ (note that both are Banach spaces), then we get $(iii)\Rightarrow(i)$ in Lemma \ref{charactprop0} and hence an equivalence there. This gives an independent proof of Theorem \ref{weightholombysequcharactsingle} (involving a dilatation parameter $c>0$). Indeed, this mentioned technique is formally more general since we only require log-convexity for $N$. But in view of Remark \ref{LCwlogrem} this difference is negligible.

\subsection{The exponential-type weight sequence setting}\label{equsectweightsequexpo}
Let $M,N$ be weight sequences such that $N_j\le AM_j$ for some $A\ge 1$ and all $j$; i.e. \eqref{strongNMrelation}. Then
$$\forall\;c>0\;\forall\;t\ge 0:\;\;\;c\omega_M(t)\le c\omega_N(t)+c\log(A)\Leftrightarrow v^c_N(t)\le A^cv^c_M(t),$$
hence
\begin{equation}\label{firstexpequ}
\forall\;c>0:\;\;\;\;H^{\infty}_{v^c_M}\subseteq H^{\infty}_{v^c_N},
\end{equation}
with continuous inclusion.

However, the natural relation in this setting is expected to be \hyperlink{ompreceqpowc}{$\preceq^{\mathfrak{c}}$}; recall Proposition \ref{firstweighfctprop}.

\begin{proposition}\label{firstexpprop}
Let $M,N\in\hyperlink{LCset}{\mathcal{LC}}$ be given. Consider the following assertions:
\begin{itemize}
\item[$(i)$] $M$ and $N$ are related by
\begin{equation}\label{assofctrelationcharlemmaequ}
\exists\;c\in\NN_{>0}\;\exists\;A\ge 1\;\forall\;j\in\NN:\;\;\;N_j\le A(M_{cj})^{1/c}.
\end{equation}

\item[$(ii)$] The weights satisfy $v_M\hyperlink{ompreceqpowc}{\preceq^{\mathfrak{c}}}v_N$.

\item[$(iii)$] We have
$$H^{\infty}_{\underline{\mathcal{M}}^{\mathfrak{c}}}\subseteq H^{\infty}_{\underline{\mathcal{N}}^{\mathfrak{c}}},\hspace{15pt}\text{and}\hspace{15pt}H^{\infty}_{\overline{\mathcal{M}}^{\mathfrak{c}}}\subseteq H^{\infty}_{\overline{\mathcal{N}}^{\mathfrak{c}}},$$
with continuous inclusions.
\end{itemize}
Then $(i)\Leftrightarrow(ii)\Rightarrow(iii)$ holds true.
\end{proposition}

\demo{Proof}
$(ii)\Rightarrow(iii)$ is immediate.

Then recall that $(v)$ in Remark \ref{LCwlogrem1} yields that $\omega_M(t)=O(\omega_N(t))$ if and only if $v_M\hyperlink{ompreceqpowc}{\preceq^{\mathfrak{c}}}v_N$. Moreover, in \cite[Lemma 6.5]{PTTvsmatrix} it is shown that $\omega_M(t)=O(\omega_N(t))$ if and only if \eqref{assofctrelationcharlemmaequ}. Combining both information $(i)\Leftrightarrow(ii)$ is verified.
\qed\enddemo

In order to proceed we have to recall and introduce some notation: For given $M\in\RR_{>0}^{\NN}$ and $c\in\NN_{>0}$ let us set (see \cite[$(2.5)$, $(2.6)$]{subaddlike} and \cite[Sect.3]{modgrowthstrange})
\begin{equation}\label{sequenceLC}
\widetilde{M}^c_j:=(M_{cj})^{1/c}.
\end{equation}
Hence $\widetilde{M}^1\equiv M$ is clear. If $M_0=1$ and $M$ is log-convex, then each $\widetilde{M}^c$ is log-convex as well, $\widetilde{M}^c_0=1$ and $\widetilde{M}^c\ge M$. If $M\in\hyperlink{LCset}{\mathcal{LC}}$, then $\widetilde{M}^c\in\hyperlink{LCset}{\mathcal{LC}}$ (for some/any $c\in\NN_{>0}$). Recall also the technical estimate \cite[Lemma 6.5, $(6.7)$]{PTTvsmatrix}:

\begin{lemma}\label{LCmoderategrowth1}
Let $M\in\hyperlink{LCset}{\mathcal{LC}}$ be given, then
\begin{equation}\label{LCmoderategrowth1equ}
\forall\;c\in\NN_{>0}\;\exists\;D\ge 1\;\forall\;t\ge 0:\;\;\;c\omega_{\widetilde{M}^c}(t)\le\omega_M(t)\le 2c\omega_{\widetilde{M}^c}(t)+Dc\Leftrightarrow\frac{1}{e^{Dc}}v^{2c}_{\widetilde{M}^c}(t)\le v_M(t)\le v^c_{\widetilde{M}^c}(t).
\end{equation}
\end{lemma}

Moreover, for any given $M\in\hyperlink{LCset}{\mathcal{LC}}$ and $c\in\NN_{>0}$ we introduce
\begin{equation}\label{auxilarysequ}
\underline{M}^c_j:=\sup_{t\ge 0}\frac{t^j}{\exp(c\omega_M(t))}=\exp(\varphi_{c\omega_M}^{*}(j)),
\end{equation}
and the second equality holds as follows (analogously as in the proof of \cite[Lemma 6.5]{PTTvsmatrix} for $\widetilde{M}^c$): Recall that by $M\in\hyperlink{LCset}{\mathcal{LC}}$ we get $\omega_M(t)=0$ for $t\in[0,1]$ (i.e. normalization) and so
\begin{align*}
\underline{M}^c_j&:=\sup_{t\ge 0}\frac{t^j}{\exp(c\omega_M(t))}=\sup_{t\ge 1}\frac{t^j}{\exp(c\omega_M(t))}=\exp(\sup_{t\ge 1}j\log(t)-c\omega_M(t))
\\&
=\exp(\sup_{s\ge 0}j s-c\omega_M(e^s))=:\exp(\varphi^{*}_{c\omega_M}(j)),
\end{align*}
where $\varphi^{*}_{\omega}$ denotes the so-called Young-conjugate of $\varphi_{\omega}: t\mapsto\omega(e^t)$. Since $\omega_M\in\hyperlink{omset0}{\mathcal{W}_0}$ we have $c\omega_M\in\hyperlink{omset0}{\mathcal{W}_0}$. By this fact and \eqref{auxilarysequ}, in particular, we know that $\underline{M}^c\in\hyperlink{LCset}{\mathcal{LC}}$. By applying \cite[Theorem 4.0.3, Lemma 5.1.3]{dissertation} (see also \cite[Lemma 2.5]{sectorialextensions}) to $\omega\equiv c\omega_M$ we have $c\omega_M\hyperlink{sim}{\sim}\omega_{\underline{M}^c}$, more precisely (by setting the weight matrix parameter $x=1$)
\begin{equation}\label{goodequivalenceclassic}
\forall\;c\in\NN_{>0}\;\exists\;D>0\;\forall\;t\ge 0:\;\;\;\omega_{\underline{M}^c}(t)\le c\omega_M(t)\le 2\omega_{\underline{M}^c}(t)+D.
\end{equation}

\begin{proposition}\label{charactprop4}
Let $M,N\in\hyperlink{LCset}{\mathcal{LC}}$ be given.
\begin{itemize}
\item[$(i)$] If we have as sets
\begin{equation}\label{charactprop4equ}
\exists\;c,d\in\NN_{>0}:\;\;\;H^{\infty}_{v^c_M}\subseteq H^{\infty}_{v^d_N},
\end{equation}
then (with the same $c$ and $c_1=2d$)
\begin{equation}\label{charactprop4equ2}
\exists\;c,c_1\in\NN_{>0}\;\exists\;C\ge 1\;\forall\;j\in\NN:\;\;\;\widetilde{N}^c_j\le C\widetilde{M}^{c_1}_j.
\end{equation}
\item[$(ii)$]
If we have the (continuous) inclusion
\begin{equation}\label{charactprop4equ1}
H^{\infty}_{\overline{\mathcal{M}}^{\mathfrak{c}}}\subseteq H^{\infty}_{\overline{\mathcal{N}}^{\mathfrak{c}}},
\end{equation}
then
\begin{equation}\label{charactprop4equ11}
\forall\;d\in\NN_{>0}\;\exists\;c\in\NN_{>0}\;\exists\;C\ge 1\;\forall\;j\in\NN:\;\;\;\widetilde{N}^d_j=(N_{jd})^{1/d}\le C(M_{jc})^{1/c}=C\widetilde{M}^c_j.
\end{equation}
\end{itemize}
\end{proposition}

\demo{Proof}
$(i)$ Let the parameters $c$ and $d$ be given and w.l.o.g. assume that $d=kc$ for some $k\in\NN_{>0}$. We apply Theorem \ref{essentialweightsequthm} to $\underline{M}^c$ and so \eqref{approximable} holds for the corresponding class. Thus there exists $f_{v_{\underline{M}^c}}\in H^{\infty}_{v_{\underline{M}^c}}$ such that $\frac{1}{v_{\underline{M}^c}}\hyperlink{sim}{\sim} t\mapsto M(f_{v_{\underline{M}^c}},t)$. By \eqref{goodequivalenceclassic} and the inclusion \eqref{charactprop4equ} one has $H^{\infty}_{v_{\underline{M}^c}}\subseteq H^{\infty}_{v^c_M}\subseteq H^{\infty}_{v^d_N}$. Consider $f^2_{v_{\underline{M}^c}}$ and so $M(f^2_{v_{\underline{M}^c}},t)=M(f_{v_{\underline{M}^c}},t)^2$. Gathering this information we arrive at
$$\exists\;A,D,C\ge 1\;\forall\;t\ge 0:\;\;\;\frac{1}{AD}\exp(c\omega_M(t))\le\frac{1}{A}\exp(2\omega_{\underline{M}^c}(t))\le M(f^2_{v_{\underline{M}^c}},t)\le C\exp(2d\omega_N(t)),$$
and for the first estimate \eqref{goodequivalenceclassic} has been used again. By applying this and \eqref{Prop32Komatsu} we obtain for all $j\in\NN$ that
\begin{align*}
M_{2dj}&=\sup_{t\ge 0}\frac{t^{2dj}}{\exp(\omega_{M}(t))}=\left(\sup_{t\ge 0}\frac{t^{2cdj}}{\exp(c\omega_{M}(t))}\right)^{1/c}\ge\frac{1}{(ADC)^{1/c}}\left(\sup_{t\ge 0}\frac{t^{2cdj}}{\exp(2d\omega_{N}(t))}\right)^{1/c}
\\&
=\frac{1}{(ADC)^{1/c}}\left(\sup_{t\ge 0}\frac{t^{cj}}{\exp(\omega_{N}(t))}\right)^{2d/c}=\frac{1}{(ADC)^{1/c}}(N_{cj})^{2d/c},
\end{align*}
and so \eqref{charactprop4equ2} is verified.\vspace{6pt}

$(ii)$ Both spaces $H^{\infty}_{\overline{\mathcal{M}}^{\mathfrak{c}}}$ and $H^{\infty}_{\overline{\mathcal{N}}^{\mathfrak{c}}}$ are Fr\'{e}chet and by the continuity of the inclusion in \eqref{charactprop4equ1} we get
$$\forall\;d>0\;\exists\;c>0\;\exists\;C\ge 1\;\forall\;f\in H^{\infty}_{\overline{\mathcal{M}}^{\mathfrak{c}}}:\;\;\; \|f\|_{v^d_N}\le C\|f\|_{v^c_M}.$$
We apply this estimate to the family of monomials $(f_k)_{k\in\NN}$ and so
$$\forall\;d>0\;\exists\;c>0\;\exists\;C\ge 1\;\forall\;k\in\NN\;\forall\;z\in\CC:\;\;\;\frac{|z|^k}{\exp(d\omega_N(|z|))}\le C\frac{|z|^k}{\exp(c\omega_M(|z|))}.$$
W.l.o.g. we can assume $c_1:=\frac{1}{c}\in\NN_{>0}, d_1:=\frac{1}{d}\in\NN_{>0}$ (with $c_1\ge d_1$) and thus, by using again \eqref{Prop32Komatsu} we obtain for all $j\in\NN$ that
\begin{align*}
N_{j/d}&=\sup_{t\ge 0}\frac{t^{j/d}}{\exp(\omega_{N}(t))}=\left(\sup_{t\ge 0}\frac{t^j}{\exp(d\omega_{N}(t))}\right)^{1/d}\le C^{1/d}\left(\sup_{t\ge 0}\frac{t^j}{\exp(c\omega_{M}(t))}\right)^{1/d}
\\&
=C^{1/d}\left(\sup_{t\ge 0}\frac{t^{j/c}}{\exp(\omega_{M}(t))}\right)^{c/d}=C^{1/d}(M_{j/c})^{c/d}.
\end{align*}
So we have verified
$$\forall\;d_1\in\NN_{>0}\;\exists\;c_1\in\NN_{>0}\;\exists\;C\ge 1\;\forall\;j\in\NN:\;\;\;(N_{jd_1})^{1/d_1}\le C(M_{jc_1})^{1/c_1},$$
i.e. \eqref{charactprop4equ11}.
\qed\enddemo

\begin{lemma}\label{tildecomparisonlemma}
Let $M,N\in\hyperlink{LCset}{\mathcal{LC}}$ be given. Then the following are equivalent:
\begin{itemize}
	\item[$(i)$] \eqref{assofctrelationcharlemmaequ} holds true,
	
	\item[$(ii$)] \eqref{charactprop4equ2} holds true,
	
	\item[$(iii)$] \eqref{charactprop4equ11} holds true.
\end{itemize}
\end{lemma}

{\itshape Note:} If in these assertions the symbol $\le$ is replaced by relation \hyperlink{preceq}{$\preceq$}, then the analogous statement holds true as well.

\demo{Proof}
Clearly, \eqref{charactprop4equ11} implies \eqref{charactprop4equ2} and \eqref{assofctrelationcharlemmaequ} (for the latter choose $d:=1$ and recall that $\widetilde{N}^1=N$).

If \eqref{charactprop4equ2} is valid, then \eqref{assofctrelationcharlemmaequ} holds true because $\widetilde{N}^c\ge N(=\widetilde{N}^1)$.

Finally, if \eqref{assofctrelationcharlemmaequ} holds, let $d\in\NN_{>0}$ be arbitrary and then $N_j\le A(M_{cj})^{1/c}$ for some $A\ge 1$ and all $j\in\NN$ implies $(N_{dj})^{1/d}\le A^{1/d}(M_{cdj})^{1/(cd)}$. This verifies \eqref{charactprop4equ11} (and hence \eqref{charactprop4equ2}).
\qed\enddemo

By combining Proposition \ref{firstexpprop}, Proposition \ref{charactprop4} and Lemma \ref{tildecomparisonlemma} we get:

\begin{theorem}\label{weightholombysequcharactpower}
Let $M,N\in\hyperlink{LCset}{\mathcal{LC}}$ be given. Then the following are equivalent:
\begin{itemize}
\item[$(i)$] $M$ and $N$ are related by \eqref{assofctrelationcharlemmaequ}, i.e.
$$\exists\;c\in\NN_{>0}\;\exists\;A\ge 1\;\forall\;j\in\NN:\;\;\;N_j\le A(M_{cj})^{1/c}=A\widetilde{M}^c_j.$$

\item[$(ii)$] The weights are related by $v_M\hyperlink{ompreceqpowc}{\preceq^{\mathfrak{c}}}v_N$.

\item[$(iii)$] We have the continuous inclusion
$$H^{\infty}_{\underline{\mathcal{M}}^{\mathfrak{c}}}\subseteq H^{\infty}_{\underline{\mathcal{N}}^{\mathfrak{c}}}.$$

\item[$(iv)$] We have the continuous inclusion
$$H^{\infty}_{\overline{\mathcal{M}}^{\mathfrak{c}}}\subseteq H^{\infty}_{\overline{\mathcal{N}}^{\mathfrak{c}}}.$$
\end{itemize}
Consequently, \hyperlink{simpowc}{$\sim^{\mathfrak{c}}$} is characterizing the equivalence of the corresponding classes (as l.c.v.s.).
\end{theorem}

\begin{remark}\label{poweroptimalrem1}
\emph{When using in the proof of $(i)$ in Proposition \ref{charactprop4} alternatively the function $\theta^c_M$ (see \eqref{charholomfctpow}), then due to the appearance of the dilatation parameter we get that the inclusion \eqref{charactprop4equ} implies (with the same $c,d$):}
\begin{equation}\label{charactprop4equ33}
	\forall\;c\in\NN_{>0}\;\exists\;d\in\NN_{>0}:\;\;\;\widetilde{N}^c\hyperlink{preceq}{\preceq}\widetilde{M}^d,
\end{equation}
\emph{equivalently}
\begin{equation}\label{charactprop4equ3}
	\exists\;c\in\NN_{>0}:\;\;\;N\hyperlink{preceq}{\preceq}\widetilde{M}^c.
\end{equation}
\emph{In order to conclude from this weaker relation we have to assume in addition that either $M$ or $N$ has \eqref{om1omegaMchar}: \eqref{charactprop4equ3} gives by definition of the associated functions
$$\exists\;c\in\NN_{>0}\;\exists\;h\ge 1\;\exists\;C\ge 0\;\forall\;t\ge 0:\;\;\;\omega_{\widetilde{M}^c}(t)\le\omega_{N}(ht)+C,$$
hence by the technical estimate \eqref{LCmoderategrowth1equ} we get
$$\frac{1}{2c}\omega_M(t)-\frac{D}{2}\le\omega_{\widetilde{M}^c}(t)\le\omega_{N}(ht)+C,$$
and so
$$\exists\;D\ge 1\;\exists\;c\in\NN_{>0}\;\exists\;h\ge 1\;\exists\;C\ge 0\;\forall\;t\ge 0:\;\;\;\omega_M(t)\le 2c\omega_{N}(ht)+2Cc+Dc.$$
\eqref{om1omegaMchar} yields \hyperlink{om1}{$(\omega_1)$} for either $\omega_M$ or $\omega_N$. An iterated application of this property implies that the previous estimate turns into $\omega_M(t)=O(\omega_N(t))$.}
\end{remark}

\subsection{The dilatation-type weight function setting}\label{equsectweightfctdila}
The strategy is to reduce this case to the weight sequence setting by applying the results from Section \ref{equsectweightsequdila} to $M^u$ (see \eqref{vBMTweight1equ1}).

\begin{proposition}\label{charactabstractweight}
Let $u$ and $w$ be normalized weight functions. Assume that $u$ is convex and that $H^{\infty}_{u}\subseteq H^{\infty}_{w}$ (as sets).
\begin{itemize}
\item[$(i)$] Then the weights are related by
\begin{equation}\label{charactabstractweightequ0}
\exists\;a\ge 1\;\forall\;t\ge 0:\;\;\;w^2(t)\le au(t),
\end{equation}
and so $u\hyperlink{ompreceqpowc}{\preceq^{\mathfrak{c}}}w$.

\item[$(ii)$] If in addition one of the weights is of moderate growth, then
\begin{equation}\label{charactabstractweightequ}
\exists\;a,b\ge 1\;\forall\;t\ge 0:\;\;\;w(at)\le bu(t),
\end{equation}
i.e. $u\hyperlink{ompreceqc}{\preceq_{\mathfrak{c}}}w$. The proof shows that we can take $a:=H$ with $H$ the constant appearing in \eqref{om6forv} for the weight being of moderate growth.
\end{itemize}
\end{proposition}

\demo{Proof}
$(i)$ By assumption and \eqref{omegavequivnewnew} we get (as sets) $H^{\infty}_{v_{M^{u}}}\subseteq H^{\infty}_{u}\subseteq H^{\infty}_{w}$. Then apply Theorem \ref{essentialweightsequthm} to $M^u$ and obtain the existence of an optimal $f_{v_{M^u}}$ such that \eqref{approximable} holds true. The previous inclusion yields now
\begin{equation}\label{charactabstractweightprotoequ}
\exists\;C,D\ge 1\;\forall\;t\ge 0:\;\;\;\omega_{M^{u}}(t)-\log(D)\le\log(C)-\log(w(t))=\log(C)+\omega^{w}(t),
\end{equation}
and recall \eqref{omegafromv} for the last equality.

There exist now two arguments to conclude; both make only use of the convexity for $u$ (see Remark \ref{convnoneed}).

1) By combining the first part of \eqref{omegavequiv} applied to $u$ and \eqref{charactabstractweightprotoequ} we get
$$\exists\;C_1\ge 1\;\exists\;A\ge 1\;\forall\;t\ge 0:\;\;\;-\log(A)+\frac{1}{2}\omega^{u}(t)\le\omega_{M^{u}}(t)\le\log(C_1)+\omega^{w}(t).$$
In view of \eqref{omegafromv} estimate \eqref{charactabstractweightequ0} is shown.

2) The second method uses again \eqref{charactabstractweightprotoequ} but follows then the proof of Proposition \ref{charactprop1}: By taking into account \eqref{vBMTweight1equ1} we get for all $j\in\NN$
\begin{align*}
M^{w}_j&=\sup_{t>0}\frac{t^j}{\exp(\omega^{w}(t))}\le CD\sup_{t>0}\frac{t^j}{\exp(\omega_{M^{u}}(t))}=CDM^{u}_j,
\end{align*}
i.e. $M^{w}$ and $M^{u}$ are related by \eqref{strongNMrelation}. Thus, by definition of associated weight functions, we have shown
$$\exists\;C_1\ge 1\;\forall\;t\ge 0:\;\;\;\omega_{M^{u}}(t)\le\omega_{M^{w}}(t)+\log(C_1)\Leftrightarrow v_{M^{w}}(t)\le C_1v_{M^{u}}(t).$$
The second estimate in \eqref{omegavequiv} applied to $w$ and the first one applied to $u$ give
\begin{equation*}\label{charactabstractweightequ0prep}
\exists\;B\ge 1\;\exists\;C_1\ge 1\;\forall\;t\ge 0:\;\;\;w(t)\le v_{M^{w}}(t)\le C_1v_{M^{u}}(t)\le BC_1\sqrt{u(t)},
\end{equation*}
which is verifying again \eqref{charactabstractweightequ0}.\vspace{6pt}

$(ii)$ If now $w$ satisfies in addition \eqref{om6forv}, then by combining this with \eqref{charactabstractweightequ0} we get
$$w(Ht)\le e^{H}w^2(t)\le e^Hau(t),$$
hence \eqref{charactabstractweightequ} is shown. If $u$ satisfies \eqref{om6forv}, then similarly $\sqrt{u(Ht)}\le e^{H/2}u(t)$ for some $H\ge 1$ and all $t\ge 0$ and so
$$w(Ht)\le\sqrt{a}\sqrt{u(Ht)}\le\sqrt{a}e^{H/2}u(t).$$
\qed\enddemo

\begin{remark}\label{charactabstractweightrem}
\emph{If in $(i)$ in Proposition \ref{charactabstractweight} we use instead of $f_{v_{M^u}}$ the function $\theta_{M^{u},1}$ from Lemma \ref{charholomfctlemma}, then a dilatation parameter $\frac{1}{2}$ appears and instead of \eqref{charactabstractweightequ0} we get the weaker estimate}
\begin{equation}\label{charactabstractweightequ0weak}
\exists\;a\ge 1\;\forall\;t\ge 0:\;\;\;(w(2t))^2\le au(t).
\end{equation}
\emph{Consequently, \eqref{charactabstractweightequ} holds then with $a:=2H$ instead of $H$. (However, in view of $(ii)$ this technical difference is negligible since there in any case a dilatation appears.)}
\end{remark}

To treat the classes $H^{\infty}_{\overline{\mathcal{V}}_{\mathfrak{c}}}$ we follow Proposition \ref{charactprop2}.

\begin{proposition}\label{charactabstractweightBeur}
Let $u$ and $w$ be normalized weights. Assume that $u$ is convex and that
\begin{equation}\label{charactabstractweightBeurequ}
H^{\infty}_{\overline{\mathcal{U}}_{\mathfrak{c}}}\subseteq H^{\infty}_{\overline{\mathcal{W}}_{\mathfrak{c}}}
\end{equation}
holds (with continuous inclusion). Then the conclusion of $(ii)$ in Proposition \ref{charactabstractweight} holds true; concerning $(i)$ the dilatation parameter $2$ in \eqref{charactabstractweightequ0weak} has to be replaced by some $a'\ge 1$ and therefore $a:=a'H$ in \eqref{charactabstractweightequ}.
\end{proposition}

\demo{Proof}
The spaces $H^{\infty}_{\overline{\mathcal{U}}_{\mathfrak{c}}}$ and $H^{\infty}_{\overline{\mathcal{W}}_{\mathfrak{c}}}$ are Fr\'{e}chet and by the continuity of the inclusion in \eqref{charactabstractweightBeurequ} we get
$$\forall\;d>0\;\exists\;c>0\;\exists\;C\ge 1\;\forall\;f\in H^{\infty}_{\overline{\mathcal{U}}_{\mathfrak{c}}}:\;\;\; \|f\|_{w_d}\le C\|f\|_{u_c}.$$
We apply this estimate again to the family of monomials $(f_k)_{k\in\NN}$ and so
$$\forall\;d>0\;\exists\;c>0\;\exists\;C\ge 1\;\forall\;k\in\NN\;\forall\;z\in\CC:\;\;\;|z|^kw(d|z|)\le C|z|^ku(c|z|).$$
Thus, by taking into account \eqref{omegafromv}, this estimate turns into
$$\forall\;d>0\;\exists\;c>0\;\exists\;C\ge 1\;\forall\;k\in\NN\;\forall\;z\in\CC:\;\;\;\frac{|z|^k}{\exp(\omega^w(d|z|))}\le C\frac{|z|^k}{\exp(\omega^u(c|z|))}.$$
Then, by using \eqref{vBMTweight1equ1}, we get for all $j\in\NN$ that
\begin{align*}
M^w_j&=\sup_{t> 0}\frac{t^j}{\exp(\omega^w(t))}=\sup_{s>0}\frac{(ds)^j}{\exp(\omega^w(ds))}\le Cd^j\sup_{s>0}\frac{s^j}{\exp(\omega^u(cs))}
\\&
=C\left(\frac{d}{c}\right)^j\sup_{s>0}\frac{s^j}{\exp(\omega^u(s))}=C\left(\frac{d}{c}\right)^jM^u_j,
\end{align*}
i.e. $M^w\hyperlink{preceq}{\preceq}M^u$ is verified. Then take $d:=1$ and the rest follows analogously as the second method in $(i)$ in Proposition \ref{charactabstractweight} and so \eqref{charactabstractweightequ0weak} holds with $a':=\frac{1}{c}$ instead of $2$ (as dilatation parameter).
\qed\enddemo

Gathering all this information we arrive at the following result:

\begin{theorem}\label{weightholombyfctcharact}
Let $u$ and $w$ be normalized weight functions. Assume that $u$ is convex and that either $u$ or $w$ is of moderate growth. Then the following are equivalent:
\begin{itemize}
\item[$(i)$] The weights satisfy $u\hyperlink{ompreceqc}{\preceq_{\mathfrak{c}}}w$.

\item[$(ii)$] We have
$$H^{\infty}_{\underline{\mathcal{U}}_{\mathfrak{c}}}\subseteq H^{\infty}_{\underline{\mathcal{W}}_{\mathfrak{c}}},$$
with continuous inclusion.

\item[$(iii)$] We have
$$H^{\infty}_{\overline{\mathcal{U}}_{\mathfrak{c}}}\subseteq H^{\infty}_{\overline{\mathcal{W}}_{\mathfrak{c}}},$$
with continuous inclusion.
\end{itemize}
If both weights are normalized and convex and at least one of them is of moderate growth, then the previous list of equivalences can be extended by
\begin{itemize}
\item[$(iv)$] The associated sequences satisfy $M^w\hyperlink{preceq}{\preceq}M^u$.

\item[$(v)$] The weights satisfy $v_{M^u}\hyperlink{ompreceqc}{\preceq_{\mathfrak{c}}}v_{M^w}$.
\end{itemize}
Consequently, for normalized convex weights of moderate growth relation \hyperlink{simc}{$\sim_{\mathfrak{c}}$} is characterizing the equivalence of the classes (as l.c.v.s.).
\end{theorem}

\demo{Proof}
$(i)\Leftrightarrow(ii)\Leftrightarrow(iii)$ follows by Propositions \ref{firstweighfctprop}, \ref{charactabstractweight} and \ref{charactabstractweightBeur}.

$(iv)\Leftrightarrow(v)$ holds by Theorem \ref{weightholombysequcharact} (even without requiring moderate growth) and $(i)\Leftrightarrow(v)$ follows by combining \eqref{omegavequivnew} with \eqref{om6forv}. Note that for this equivalence we have to use convexity for {\itshape both} weights whereas it is enough to assume moderate growth for one of the weights; see the proof of $(ii)$ in Proposition \ref{charactabstractweight}.
\qed\enddemo

\begin{remark}\label{Clunieremark}
\emph{We have reduced the proof of Theorem \ref{weightholombyfctcharact} to the weight sequence setting by involving the associated weight sequence $M^u$. For this procedure the appearing technical growth assumption \eqref{om6forv} seems to be unavoidable in order to get control on the quadratic power in \eqref{omegavequiv}. Similarly, in the ultradifferentiable setting moderate growth is needed when transferring techniques from the weight function to the weight sequence setting.}\vspace{6pt}

\emph{In \cite{integralprescribed} the author has also started with an abstractly given increasing function $\phi$ such that $\phi(t)\neq O(\log(t))$ and that $t\mapsto\phi(e^t)$ is convex. Note that the first requirement is weaker than $\phi(t)=o(\log(t))$ which always holds true for $\omega^v$ when $v$ is a weight function, and the second basic assumption holds true for $\omega^v$ if and only if $v$ is convex (see \eqref{vconvexity}). Summarizing, we have the relation $\phi\equiv-\log\circ v$, with $v$ being a weight function.}

\emph{Moreover, the representation $\phi(t)=\int_1^t\frac{\psi(s)}{s}ds$ is used which should be compared with \eqref{assointrepr}; i.e. $\psi\equiv\Sigma_M$ with $\Sigma_M$ denoting the counting function in \eqref{counting}. By using $\psi$ the sequence $r=(r_i)_{i\ge 1}$ is defined via $\psi(i)=r_i$ and finally $a_n:=\frac{1}{r_1\cdots r_n}$, i.e. $r$ corresponds to the quotient sequence $\mu^v$ and $a_n=(M^v_n)^{-1}$. Note that $\psi(1)=0$ precisely corresponds to normalization of given $v$.}

\emph{Then, in \cite[Sect. 2]{integralprescribed} it is shown that there exists an entire function $f$ such that
$$\lim_{t\rightarrow+\infty}\frac{\log(M(f,t))}{\phi(t)}=1.$$
In fact, in \cite[Lemma 2]{integralprescribed} it is shown that $\widetilde{f}:=f|_{\RR}$ satisfies}
\begin{equation}\label{Clunieremarkequ}
\lim_{t\rightarrow+\infty}\frac{\log(M(f,t))}{\phi(t)}=\lim_{t\rightarrow+\infty}\frac{\log(\widetilde{f}(t))}{\phi(t)}=1.
\end{equation}
\emph{\eqref{Clunieremarkequ} means that for all $1>\epsilon>0$ we have $\widetilde{f}(t)\ge\exp((1-\epsilon)\phi(t))$ for all $t\ge 0$ large enough and so the inclusion $H^{\infty}_{v}\subseteq H^{\infty}_{w}$ implies (by the correspondence $\phi_v\equiv-\log\circ v$, $\phi_w\equiv-\log\circ w$)}
$$\exists\;C\ge 1\;\forall\;\epsilon>0\;\exists\;D\ge 1\;\forall\;t\ge 0:\;\;\;\frac{1}{D}(v(t))^{\epsilon-1}=\frac{1}{D}\exp((1-\epsilon)\phi_v(t))\le C\exp(\phi_w(t))=C\frac{1}{w(t)}.$$
\emph{But, since $1-\epsilon<1$, in order to verify \eqref{charactabstractweightequ}, property \eqref{om6forv} has to be applied again (to one of the weights).}
\end{remark}

\subsection{The exponential-type weight function setting}\label{equsectweightfctpow}
Recall that by Proposition \ref{firstweighfctprop} for given weights $v$ and $w$ being related via $v\hyperlink{ompreceqpowc}{\preceq^{\mathfrak{c}}}w$ we get
$$H^{\infty}_{\underline{\mathcal{V}}^{\mathfrak{c}}}\subseteq H^{\infty}_{\underline{\mathcal{W}}^{\mathfrak{c}}},\hspace{15pt}H^{\infty}_{\overline{\mathcal{V}}^{\mathfrak{c}}}\subseteq H^{\infty}_{\overline{\mathcal{W}}^{\mathfrak{c}}},$$
with continuous inclusions.

\begin{proposition}\label{charactprop4forweight}
Let $u$ and $w$ be normalized weights. Assume that $u$ is convex and that either
\begin{itemize}
\item[$(i)$] we have (as sets)
\begin{equation}\label{charactprop4forweightequ}
\exists\;c,d\in\NN_{>0}:\;\;\;H^{\infty}_{u^c}\subseteq H^{\infty}_{w^d},
\end{equation}
or that
\item[$(ii)$]
\begin{equation}\label{charactprop4forweightequ1}
H^{\infty}_{\overline{\mathcal{U}}^{\mathfrak{c}}}\subseteq H^{\infty}_{\overline{\mathcal{W}}^{\mathfrak{c}}}
\end{equation}
holds (with continuous inclusion).
\end{itemize}
Then
\begin{equation}\label{charactprop4forweightequ2}
\exists\;a\ge 1\;\exists\;b\ge 1\;\forall\;t\ge 0:\;\;\;w^a(t)\le bu(t)
\end{equation}
is valid, i.e. $u\hyperlink{ompreceqpowc}{\preceq^{\mathfrak{c}}}w$.
\end{proposition}

\demo{Proof}
$(i)$ We follow the proof of $(i)$ in Proposition \ref{charactprop4}. Let the parameters $c$ and $d$ be given (w.l.o.g. $d\ge c\ge 1$). Then apply Theorem \ref{essentialweightsequthm} to $\underline{M^u}^c$ (recall \eqref{auxilarysequ}) and so \eqref{approximable} holds for the corresponding class. Thus there exists $f_{v_{\underline{M^u}^c}}\in H^{\infty}_{v_{\underline{M^u}^c}}$ such that $\frac{1}{v_{\underline{M^u}^c}}\hyperlink{sim}{\sim} t\mapsto M(f_{v_{\underline{M^u}^c}},t)$. The inclusion \eqref{charactprop4forweightequ} together with \eqref{omegavequivnewnew} and \eqref{goodequivalenceclassic} give $H^{\infty}_{v_{\underline{M^u}^c}}\subseteq H^{\infty}_{v^c_{M^u}}\subseteq H^{\infty}_{u^c}\subseteq H^{\infty}_{w^d}$. Consider the function $f^2_{v_{\underline{M^u}^c}}$ and so gathering all this information we arrive at
$$\exists\;A,D,C\ge 1\;\forall\;t\ge 0:\;\;\;\frac{1}{AD}\exp(c\omega_{M^u}(t))\le\frac{1}{A}\exp(2\omega_{\underline{M^u}^c}(t))\le M(f^2_{v_{\underline{M^u}^c}},t)\le C(w^d(t))^{-2},$$
where for the first estimate \eqref{goodequivalenceclassic} applied to $M^u$ has been used again. Thus, by the first part of \eqref{omegavequiv} applied to $u$ and for which {\itshape convexity} is needed, we have
$$\exists\;A,B,C,D\ge 1\;\forall\;t\ge 0:\;\;\;\frac{1}{ADB^c}(u^{c/2}(t))^{-1}\le\frac{1}{AD}\exp(c\omega_{M^u}(t))\le C(w^d(t))^{-2},$$
hence $(w^{4d/c}(t))\le(AB^cCD)^{2/c}u(t)$ follows; i.e. \eqref{charactprop4forweightequ2}.

Alternatively, the above estimate yields
\begin{align*}
&\exists\;A,C,D\ge 1\;\forall\;t\ge 0:
\\&
c\omega_{M^u}(t)-\log(AD)\le\log(C)-2d\log(w(t))=\log(C)+2d\omega^w(t),
\end{align*}
and then in view of \eqref{vBMTweight1equ1} we can follow part $(i)$ in Proposition \ref{charactprop4} and get $\widetilde{M^w}^c_j\le C_1\widetilde{M^u}^{2d}_j$; i.e. \eqref{charactprop4equ2} for $M^w$ and $M^u$. By Lemma \ref{tildecomparisonlemma} the sequences $M^w$ and $M^u$ are related by \eqref{assofctrelationcharlemmaequ} and so $(i)\Leftrightarrow(ii)$ in Proposition \ref{firstexpprop} yields $v_{M^u}\hyperlink{ompreceqpowc}{\preceq^{\mathfrak{c}}}v_{M^w}$. Finally \eqref{omegavequivnew} implies \eqref{charactprop4forweightequ2}. Note that for this last step again only convexity for $u$ is required.\vspace{6pt}

$(ii)$ Both spaces $H^{\infty}_{\overline{\mathcal{U}}^{\mathfrak{c}}}$ and $H^{\infty}_{\overline{\mathcal{W}}^{\mathfrak{c}}}$ are Fr\'{e}chet and by the continuity of the inclusion in \eqref{charactprop4forweightequ1} we get
$$\forall\;d>0\;\exists\;c>0\;\exists\;C\ge 1\;\forall\;f\in H^{\infty}_{\overline{\mathcal{U}}^{\mathfrak{c}}}:\;\;\; \|f\|_{w^d}\le C\|f\|_{u^c}.$$
We apply this estimate to the family of monomials $(f_k)_{k\in\NN}$ and so
\begin{align*}
&\forall\;d>0\;\exists\;c>0\;\exists\;C\ge 1\;\forall\;k\in\NN\;\forall\;z\in\CC:
\\&
\frac{|z|^k}{\exp(d\omega^w(|z|))}=|z|^kw^d(|z|)\le C|z|^ku^c(|z|)=C\frac{|z|^k}{\exp(c\omega^u(|z|))}.
\end{align*}
W.l.o.g. we can assume $c_1:=\frac{1}{c}\in\NN_{>0}, d_1:=\frac{1}{d}\in\NN_{>0}$ (with $c_1\ge d_1$) and thus, by taking into account \eqref{vBMTweight1equ1} we follow the proof of $(ii)$ in Proposition \ref{charactprop4} and get \eqref{charactprop4equ11} between $M^w$ and $M^u$. The rest follows now, by involving again Lemma \ref{tildecomparisonlemma}, as in the second method in case $(i)$.
\qed\enddemo

Summarizing all information we obtain the following characterization:

\begin{theorem}\label{weightholomcharactpower}
Let $u$ and $w$ be normalized weights and assume that $u$ is convex. Then the following are equivalent:
\begin{itemize}
\item[$(i)$] The weights satisfy $u\hyperlink{ompreceqc}{\preceq^{\mathfrak{c}}}w$.

\item[$(ii)$] We have
$$H^{\infty}_{\underline{\mathcal{U}}^{\mathfrak{c}}}\subseteq H^{\infty}_{\underline{\mathcal{W}}^{\mathfrak{c}}},$$
with continuous inclusion.

\item[$(iii)$] We have
$$H^{\infty}_{\overline{\mathcal{U}}^{\mathfrak{c}}}\subseteq H^{\infty}_{\overline{\mathcal{W}}^{\mathfrak{c}}},$$
with continuous inclusion.
\end{itemize}
Consequently, for convex weights relation \hyperlink{simpowc}{$\sim^{\mathfrak{c}}$} is characterizing the equivalence of the corresponding classes (as l.c.v.s.).
Moreover, if both $u$ and $w$ are convex, then the previous list of equivalences can be extended by
\begin{itemize}
\item[$(iv)$] $M^u$ and $M^w$ are related by \eqref{assofctrelationcharlemmaequ}; i.e.
$$\exists\;c\in\NN_{>0}\;\exists\;A\ge 1\;\forall\;j\in\NN:\;\;\;M^w_j\le A(M^u_{cj})^{1/c}.$$

\item[$(v)$] The weights are related by $v_{M^u}\hyperlink{ompreceqpowc}{\preceq^{\mathfrak{c}}}v_{M^w}$.
\end{itemize}
\end{theorem}

\demo{Proof}
$(i)\Leftrightarrow(ii)\Leftrightarrow(iii)$ is valid by Propositions \ref{firstweighfctprop} and \ref{charactprop4forweight}. Moreover, $(iv)\Leftrightarrow(v)$ is precisely $(i)\Leftrightarrow(ii)$ in Proposition \ref{firstexpprop} applied to $M^u$ and $M^w$. Finally, $(i)\Leftrightarrow(v)$ holds by \eqref{omegavequivnew}; i.e. $v_{M^u}\hyperlink{simpowc}{\sim^{\mathfrak{c}}}u$, $v_{M^w}\hyperlink{simpowc}{\sim^{\mathfrak{c}}}w$. For this convexity is required for {\itshape both} weights; more precisely convexity for $u$ is needed for $(v)\Rightarrow(i)$ and for $w$ concerning $(i)\Rightarrow(v)$.
\qed\enddemo

We close with the following observation:

\begin{remark}\label{expweightfctfinalrem}
\emph{When using in the proof of $(i)$ in Proposition \ref{charactprop4forweight} the function $\theta^c_{M^u}$ from \eqref{charholomfctpow}, then the appearance of a dilatation parameter requires to involve \eqref{om1forv} (for one of the weights) in order to obtain the conclusion. For this recall Remark \ref{poweroptimalrem1} and note that \eqref{om1forv} for either $u$ or $w$ is equivalent to having \eqref{om1omegaMchar} for either $M^u$ or $M^w$; see Lemma \ref{om1rem}.}
\end{remark}

\section{Stability under point-wise multiplication}\label{pointwisesection}
The aim of this section is to characterize the stability of weighted classes of entire functions under the point-wise multiplication of functions. More precisely, let us consider the (bilinear) multiplication operator
$$\mathfrak{m}:\;\;\;(f,g)\mapsto f\cdot g,$$
acting on weighted spaces of entire functions.

\subsection{Preliminaries}
We start with the following observation:

\begin{proposition}\label{poinwisepowerclear}
Let $v$ be a (normalized) weight function. Then both classes $H^{\infty}_{\underline{\mathcal{V}}^{\mathfrak{c}}}$ and $H^{\infty}_{\overline{\mathcal{V}}^{\mathfrak{c}}}$ are automatically closed under the action of $\mathfrak{m}$ and $$\mathfrak{m}:\;H^{\infty}_{\underline{\mathcal{V}}^{\mathfrak{c}}}\times H^{\infty}_{\underline{\mathcal{V}}^{\mathfrak{c}}}\rightarrow H^{\infty}_{\underline{\mathcal{V}}^{\mathfrak{c}}},\;\;\;\text{resp.}\;\;\;\mathfrak{m}:\;H^{\infty}_{\overline{\mathcal{V}}^{\mathfrak{c}}}\times H^{\infty}_{\overline{\mathcal{V}}^{\mathfrak{c}}}\rightarrow H^{\infty}_{\overline{\mathcal{V}}^{\mathfrak{c}}},$$
is continuous.
\end{proposition}

\demo{Proof}
Let $f\in H^{\infty}_{v^c}$, $g\in H^{\infty}_{v^d}$ for some $c,d>0$. Then $f\cdot g\in H^{\infty}_{v^{c+d}}$ holds true. The continuity of $\mathfrak{m}$ follows by \cite[\S 40, Sect. 2, $(1)$, $(11)$]{Koethe79} and \cite[\S 40, Sect. 5, $(3)$, $(4)$]{Koethe79}.
\qed\enddemo

On the other hand we have the following statement which has been sent to the author in a private communication by Prof. Jos\'{e} Bonet.

\begin{proposition}\label{Pepecounterproof}
Let $v$ be an essential weight. Then $H^{\infty}_v$ is not closed under point-wise multiplication.
\end{proposition}

\demo{Proof}
Suppose that $\mathfrak{m}: H^{\infty}_v\rightarrow H^{\infty}_v$ is continuous, so
$$\exists\;C>0\;\forall\;f\in H^\infty_v:\;\;\;\|f^2\|_v \le C\|f\|_v.$$
Since $v$ is essential we have the following: Given an arbitrary $x\in\CC$ we find $g\in H^{\infty}_v$ such that $|g(z)|\le\frac{1}{v(z)}$ for all $z\in\CC$ and $g(x)=\frac{1}{v(x)}$, i.e. $\|g\|_v=1$.

Then
$$\frac{1}{v(x)}\le v(x)|g(x)^2|\le\|g^2\|_v \le C\|g\|_v = C,$$
and because $x$ is arbitrary we conclude $\frac{1}{C}\leq v(x)$ for all $x\in\CC$.

Let $f\in H^\infty_v$ be arbitrary and so there exists $M>0$ such that $v(z)|f(z)|\le M$ for all $z\in\CC$. By the previous estimate $|f(z)|\le MC$ for all $z\in\CC$ follows and since $f$ is entire this boundedness implies that $f$ is constant, a contradiction. (Recall that, since $v$ is rapidly decreasing, at least all polynomials are contained in $H^\infty_v$.)
\qed\enddemo

The last argument in the proof before shall be compared with the characterization \cite[Prop. 2.1]{BonetDomanskiLindstroem99} for the unit disc case.

\subsection{Convolved weight sequences}\label{convolvedsect}
We are gathering now some technical comments which are needed later on. Let $M,N\in\RR_{>0}^{\NN}$, then define the {\itshape convolved sequence} $M\star N$ by
\begin{equation}\label{convolvesequ}
M\star N_j:=\min_{0\le k\le j}M_kN_{j-k},\;\,\;j\in\NN,
\end{equation}
see \cite[$(3.15)$]{Komatsu73}. Hence, obviously $M\star N=N\star M$ and $(M\star N)_0=1$ provided that $M_0=N_0=1$.

\begin{remark}\label{convolvesequrem0}
\emph{In \cite[Lemma 3.5]{Komatsu73} the following is shown for log-convex weight sequences $M,N$ (resp. $M,N\in\hyperlink{LCset}{\mathcal{LC}}$):}

\begin{itemize}
\item[$(*)$] \emph{$M\star N$ is also a log-convex weight sequence (resp. $M\star N\in\hyperlink{LCset}{\mathcal{LC}}$).}

\item[$(*)$] \emph{The corresponding quotient sequence $\mu\star\nu$ is obtained when rearranging resp. ordering the sequences $\mu$ and $\nu$ in the order of growth.}

\item[$(*)$] \emph{This identity yields by definition}
$$\forall\;t\ge 0:\;\;\;\Sigma_{M\star N}(t)=\Sigma_M(t)+\Sigma_N(t),$$
\emph{and so by \eqref{assointrepr}}
$$\forall\;t\ge 0:\;\;\;\omega_{M\star N}(t)=\omega_M(t)+\omega_N(t).$$

\item[$(*)$] \emph{Thus, by recalling \eqref{weights}, we have}
\begin{equation}\label{weightsconv}
\forall\;c>0\;\forall\;t\ge 0:\;\;\;v_{M\star N,c}(t)=\exp(-\omega_{M\star N}(ct))=\exp(-\omega_M(ct)-\omega_N(ct))=v_{M,c}(t)v_{N,c}(t).
\end{equation}
\end{itemize}
\end{remark}

\begin{remark}\label{convolvesequrem}
\emph{Let $M,N\in\RR_{>0}^{\NN}$, we summarize more properties for $M\star N$.}

\begin{itemize}
\item[$(i)$] \emph{For any $c>0$ by definition we have $M^c\star N^c=(M\star N)^c$; see $(a)$ in Remark \ref{powermultisequ}.}

\item[$(ii)$] \emph{For all $j\in\NN$ we have $(M\star N)_j\le\min\{M_0N_j,N_0M_j\}$. So, if in addition $M_0=N_0=1$, then we get $M\star N\le\min\{M,N\}$.}

    \emph{On the other hand $M\hyperlink{preceq}{\preceq}M\star M$ amounts to requiring \hyperlink{mg}{$(\on{mg})$} (and $N\hyperlink{preceq}{\preceq}M\star M$ yields the mixed version of \hyperlink{mg}{$(\on{mg})$} between $N$ and $M$).}

    \emph{Summarizing, for any $M\in\RR_{>0}^{\NN}$ we have that $M$ and $M\star M$ are equivalent if and only if $M$ satisfies \hyperlink{mg}{$(\on{mg})$}.}

\item[$(iii)$] \emph{If $M$ is a log-convex weight sequence, then}
\begin{equation}\label{preom6}
\forall\;t\ge 0:\;\;\;\omega_M(t)\le 2\omega_M(t)=\omega_{M\star M}(t),
\end{equation}
\emph{i.e.
$$\forall\;c>0\;\forall\;t\ge 0:\;\;\;v_{M,c}(t)^2=v_{M\star M,c}(t)\le v_{M,c}(t).$$
And $M$ satisfies \hyperlink{mg}{$(\on{mg})$}, equivalently $M\hyperlink{preceq}{\preceq}M\star M$, if and only if}
\begin{equation}\label{om6asso}
\exists\;H\ge 1\;\forall\;t\ge 0:\;\;\;\omega_{M\star M}(t)=2\omega_M(t)\le\omega_M(Ht)+H,
\end{equation}
\emph{i.e. \hyperlink{om6}{$(\omega_6)$} for $\omega_M$, see \cite[Prop. 3.6]{Komatsu73}. In this situation, \eqref{weightsconv} turns into}
\begin{equation}\label{weightsconvmg}
\exists\;H\ge 1\;\forall\;c>0\;\forall\;t\ge 0:\;\;\;v_{M,cH}(t)\le e^H v_{M\star M,c}(t)=e^Hv_{M,c}(t)^2\le e^Hv_{M,c}(t).
\end{equation}
\end{itemize}
\end{remark}

\subsection{The dilatation-type weight sequence setting}\label{multiclosedweightsequ}
Let $M,N$ be log-convex weight sequences and $f\in H^{\infty}_{v_{M,c_1}}$, $g\in H^{\infty}_{v_{N,c_2}}$. Then by Remark \ref{convolvesequrem0} we get
\begin{align*}
\sup_{z\in\CC}|(f\cdot g)(z)|&=\sup_{z\in\CC}|f(z)|\cdot|g(z)|\le AB\exp(\omega_M(c_1|z|)+\omega_N(c_2|z|))
\\&
\le C\exp(\omega_M(d|z|)+\omega_N(d|z|))=C\exp(\omega_{M\star N}(d|z|)),
\end{align*}
when taking $C:=AB$ and $d:=\max\{c_1,c_2\}$. Thus,
\begin{equation}\label{multiplieroperator}
\mathfrak{m}: H^{\infty}_{v_{M,c_1}}\times H^{\infty}_{v_{N,c_2}}\longrightarrow H^{\infty}_{v_{M\star N,\max\{c_1,c_2\}}},\hspace{15pt}(f,g)\mapsto f\cdot g,
\end{equation}
and these observations yield the following result:

\begin{proposition}\label{convolutorlemma1}
Let $M$ be a log-convex weight sequence with \hyperlink{mg}{$(\on{mg})$}. Then
$$\exists\;H\ge 1\;\forall\;c_1,c_2>0\;\forall\;f\in H^{\infty}_{v_{M,c_1}},g\in H^{\infty}_{v_{M,c_2}}:\;\;\;f\cdot g\in H^{\infty}_{v_{M,H\max\{c_1,c_2\}}},$$
with $H$ being the constant appearing in \hyperlink{om6}{$(\omega_6)$} for $\omega_M$.

Consequently, both $H^{\infty}_{\underline{\mathcal{M}}_{\mathfrak{c}}}$ and $H^{\infty}_{\overline{\mathcal{M}}_{\mathfrak{c}}}$ are closed under point-wise multiplication and the multiplication operator $\mathfrak{m}$ is continuously acting on these spaces.
\end{proposition}

\demo{Proof}
In \eqref{multiplieroperator} we take $M=N$ and by \eqref{weightsconvmg} we are done. The continuity for $\mathfrak{m}$ follows again by \cite[\S 40, Sect. 2, $(1)$, $(11)$]{Koethe79} and \cite[\S 40, Sect. 5, $(3)$, $(4)$]{Koethe79}.
\qed\enddemo

The next result deals with the converse implication.

\begin{proposition}\label{convolutorlemma2}
Let $M,N,L$ be log-convex weight sequences.
\begin{itemize}
\item[$(i)$] Assume that
$$\exists\;c_1,c_2,d>0\;\forall\;f\in H^{\infty}_{v_{M,c_1}},g\in H^{\infty}_{v_{N,c_2}}:\;\;\;f\cdot g\in H^{\infty}_{v_{L,d}}.$$
Then the sequences satisfy
$$L\hyperlink{preceq}{\preceq}M\star N.$$

\item[$(ii)$] If $(i)$ holds for $M\equiv N\equiv L$, then $M$ has to satisfy \hyperlink{mg}{$(\on{mg})$}.
\end{itemize}
Alternatively, in order to conclude in the assumption in $(i)$ we can replace some or all $\exists$ by $\forall$.
\end{proposition}

\demo{Proof}
$(i)$ Let $c_1,c_2,d>0$ be given and by assumption for any $f\in H^{\infty}_{v_{M,c_1}}$, $g\in H^{\infty}_{v_{N,c_2}}$ we get
$$\exists\;C\ge 1\;\forall\;z\in\CC:\;\;\;|f(z)\cdot g(z)|\le C\exp(\omega_{L}(d|z|)).$$
We apply this estimate to $f\equiv \theta_{M,c_1}$ and $g\equiv\theta_{N,c_2}$ from Lemma \ref{charholomfctlemma} and so
$$\forall\;t\ge 0:\;\;\;\exp(\omega_M(c_1t/2))\exp(\omega_N(c_2t/2))\le|\theta_{M,c_1}(t)\cdot\theta_{N,c_2}(t)|\le C\exp(\omega_L(dt)).$$
Take $c:=\min\{c_1,c_2\}$ and then Remark \ref{convolvesequrem0} gives
$$\omega_{M\star N}(ct/2)=\omega_M(ct/2)+\omega_N(ct/2)\le\omega_M(c_1t/2)+\omega_N(c_2t/2)\le\log(C)+\omega_L(dt).$$
Since $M\star N$ is also a log-convex weight sequence (see again Remark \ref{convolvesequrem0}) we can apply \eqref{Prop32Komatsu} and get for all $j\in\NN$
\begin{align*}
L_j&=\sup_{t\ge 0}\frac{t^j}{\exp(\omega_{L}(t))}=\sup_{s\ge 0}\frac{(ds)^j}{\exp(\omega_{L}(ds))}\le Cd^j\sup_{s\ge 0}\frac{s^j}{\exp(\omega_{M\star N}(cs/2))}
\\&
=C\left(\frac{2d}{c}\right)^j\sup_{u\ge 0}\frac{u^j}{\exp(\omega_{M\star N}(u))}=C\left(\frac{2d}{c}\right)^j(M\star N)_j,
\end{align*}
i.e. $L\hyperlink{preceq}{\preceq}M\star N$.\vspace{6pt}

$(ii)$ If we take $M\equiv N\equiv L$ in $(i)$, then $M\hyperlink{preceq}{\preceq}M\star M$ and this precisely means \hyperlink{mg}{$(\on{mg})$} for $M$; see $(ii)$ in Remark \ref{convolvesequrem}.
\qed\enddemo

For the class $H^{\infty}_{\overline{\mathcal{M}}_{\mathfrak{c}}}$ we cannot make use of the optimal functions $\theta_{M,c}$.

\begin{proposition}\label{convolutorlemma3}
Let $M$ be a log-convex weight sequence. Assume that $H^{\infty}_{\overline{\mathcal{M}}_{\mathfrak{c}}}$ is closed under point-wise multiplication and that $\mathfrak{m}$ is continuously acting on this space. Then $M$ has to satisfy \hyperlink{mg}{$(\on{mg})$}.
\end{proposition}

\demo{Proof}
The proof is inspired by \cite[Prop. 4.7]{Borelmapalgebraity}. By assumption
$$\mathfrak{m}_{\Delta}:\;\;\;H^{\infty}_{\overline{\mathcal{M}}_{\mathfrak{c}}}\rightarrow H^{\infty}_{\overline{\mathcal{M}}_{\mathfrak{c}}},\hspace{15pt}f\mapsto f^2,$$
is continuous. This amounts to $H^{\infty}_{\overline{\mathcal{M}}_{\mathfrak{c}}}\subseteq H^{2,\infty}_{\overline{\mathcal{M}}_{\mathfrak{c}}}$ for the Fr\'{e}chet space
$$H^{2,\infty}_{\overline{\mathcal{M}}_{\mathfrak{c}}}:=\{f\in H(\CC): \|f^2\|_{v_{M,c}}<+\infty,\;\forall\; c>0\}.$$
By the closed graph theorem the inclusion $H^{\infty}_{\overline{\mathcal{M}}_{\mathfrak{c}}}\subseteq H^{2,\infty}_{\overline{\mathcal{M}}_{\mathfrak{c}}}$ is continuous and hence
\begin{equation}\label{convolutorlemma3equ}
\forall\;d>0\;\exists\;c>0\;\exists\;C\ge 1\;\forall\;f\in H^{\infty}_{\overline{\mathcal{M}}_{\mathfrak{c}}}:\;\;\;\|f^2\|_{v_{M,d}}\le C\|f\|_{v_{M,c}}.
\end{equation}
We apply this estimate to the monomials $f_k=z^k$ and get:
$$\forall\;d>0\;\exists\;c>0\;\exists\;C\ge 1\;\forall\;k\in\NN\;\forall\;z\in\CC:\;\;\;\frac{|z|^{2k}}{\exp(\omega_M(d|z|))}\le C\frac{|z|^k}{\exp(\omega_M(c|z|))}.$$
\eqref{Prop32Komatsu} gives for all $j\in\NN$ that
\begin{align*}
M_{2j}&=\sup_{t\ge 0}\frac{t^{2j}}{\exp(\omega_{M}(t))}=\sup_{s\ge 0}\frac{(ds)^{2j}}{\exp(\omega_{M}(ds))}\le Cd^{2j}\sup_{s\ge 0}\frac{s^j}{\exp(\omega_{M}(cs))}
\\&
=C\left(\frac{d^2}{c}\right)^j\sup_{u\ge 0}\frac{u^j}{\exp(\omega_{M}(u))}=C\left(\frac{d^2}{c}\right)^jM_j.
\end{align*}
So far we have verified
\begin{equation}\label{mgdiag}
\exists\;A,B\ge 1\;\forall\;j\in\NN:\;\;\;M_{2j}\le AB^jM_j.
\end{equation}
Hence, by \cite[Thm. 9.5.1]{dissertation} applied to the constant matrix $\{M\}$ (see also \cite[Thm. 1]{matsumoto}), \eqref{mgdiag} is equivalent to the fact that $M$ has to satisfy \hyperlink{mg}{$(\on{mg})$}.
\qed\enddemo

By combining Propositions \ref{convolutorlemma1}, \ref{convolutorlemma2} and \ref{convolutorlemma3} we arrive at the following characterization:

\begin{theorem}\label{convolutorthm}
Let $M$ be a log-convex weight sequence. Then the following are equivalent:
\begin{itemize}
\item[$(i)$] $H^{\infty}_{\underline{\mathcal{M}}_{\mathfrak{c}}}$ is closed under point-wise multiplication and $\mathfrak{m}$ is continuously acting on this space.

\item[$(ii)$] $H^{\infty}_{\overline{\mathcal{M}}_{\mathfrak{c}}}$ is closed under point-wise multiplication and $\mathfrak{m}$ is continuously acting on this space.

\item[$(iii)$] $M$ has \hyperlink{mg}{$(\on{mg})$}.
\end{itemize}
\end{theorem}

We close this section with the following observation.

\begin{remark}\label{mixedpointwise}
\emph{If $M$ is violating \hyperlink{mg}{$(\on{mg})$}, then one is able to infer a mixed version of Theorem \ref{convolutorthm}: When taking $M=N$ in \eqref{multiplieroperator} we see that both $\mathfrak{m}: H^{\infty}_{\underline{\mathcal{M}}_{\mathfrak{c}}}\times H^{\infty}_{\underline{\mathcal{M}}_{\mathfrak{c}}} \rightarrow H^{\infty}_{\underline{\mathcal{L}}_{\mathfrak{c}}}$ and $\mathfrak{m}: H^{\infty}_{\overline{\mathcal{M}}_{\mathfrak{c}}}\times H^{\infty}_{\overline{\mathcal{M}}_{\mathfrak{c}}}\rightarrow H^{\infty}_{\overline{\mathcal{L}}_{\mathfrak{c}}}$ are continuous with $L:=M\star M$.}

\emph{On the other hand, the choice $M=N$ in $(i)$ in Proposition \ref{convolutorlemma2} gives $L\hyperlink{preceq}{\preceq}M\star M$; i.e. the mixed version of \hyperlink{mg}{$(\on{mg})$} between $M$ and $L$. Similarly, when assuming that $\mathfrak{m}: H^{\infty}_{\overline{\mathcal{M}}_{\mathfrak{c}}}\times  H^{\infty}_{\overline{\mathcal{M}}_{\mathfrak{c}}}\rightarrow H^{\infty}_{\overline{\mathcal{L}}_{\mathfrak{c}}}$ is continuous, then the proof of Proposition \ref{convolutorlemma3} gives the mixed version of \eqref{mgdiag} and \cite[Thm. 9.5.1]{dissertation} implies $L\hyperlink{preceq}{\preceq}M\star M$.}

\emph{In both cases $L\hyperlink{preceq}{\preceq}M\star M$ yields $H^{\infty}_{\underline{\mathcal{M}\star\mathcal{M}}_{\mathfrak{c}}}\subseteq H^{\infty}_{\underline{\mathcal{L}}_{\mathfrak{c}}}$ and $H^{\infty}_{\overline{\mathcal{M}\star\mathcal{M}}_{\mathfrak{c}}}\subseteq H^{\infty}_{\overline{\mathcal{L}}_{\mathfrak{c}}}$; recall Proposition \ref{charactprop}.}
\end{remark}

\subsection{The dilatation-type weight function setting}\label{multiclosedweightfct}
Let $u$ and $v$ be weight functions and $f\in H^{\infty}_{u_{c_1}}$, $g\in H^{\infty}_{v_{c_2}}$. Then we get $f\cdot g\in H^{\infty}_{w_d}$ with $d:=\max\{c_1,c_2\}$ and $w:=u\cdot v$ by the same estimate as in Section \ref{multiclosedweightsequ} above. Thus, the bilinear operator
$$\mathfrak{m}: H^{\infty}_{u_{c_1}}\times H^{\infty}_{v_{c_2}}\longrightarrow H^{\infty}_{w_{\max\{c_1,c_2\}}},\hspace{15pt}(f,g)\mapsto f\cdot g,$$
is well-defined and continuous.

When taking $u\equiv v$, then we get the following result which is analogous to Proposition \ref{convolutorlemma1}.

\begin{proposition}\label{convolutorlemmaweight1}
Let $v$ be a (normalized) weight function of moderate growth. Then
$$\exists\;H\ge 1\;\forall\;c_1,c_2>0\;\forall\;f\in H^{\infty}_{v_{c_1}}, g\in H^{\infty}_{v_{c_2}}:\;\;\;f\cdot g\in H^{\infty}_{v_{H\max\{c_1,c_2\}}},$$
and $H$ can be chosen to be the constant appearing in \eqref{om6forv}. Consequently, both $H^{\infty}_{\underline{\mathcal{V}}_{\mathfrak{c}}}$ and $H^{\infty}_{\overline{\mathcal{V}}_{\mathfrak{c}}}$ are closed under point-wise multiplication and the multiplication operator $\mathfrak{m}$ is continuously acting on these spaces.
\end{proposition}

\demo{Proof}
Let $f\in H^{\infty}_{v_{c_1}}$, $g\in H^{\infty}_{v_{c_2}}$ be given and so $f\cdot g\in H^{\infty}_{v_c\cdot v_c}$ with $c:=\max\{c_1,c_2\}$. By \eqref{om6forv} we find $v(Ht)\le e^Hv(t)^2$ for some $H\ge 1$ and all $t\ge 0$, hence $f\cdot g\in H^{\infty}_{v_{Hc}}$.

The continuity for $\mathfrak{m}$ follows again by \cite[\S 40, Sect. 2, $(1)$, $(11)$]{Koethe79} and \cite[\S 40, Sect. 5, $(3)$, $(4)$]{Koethe79}.
\qed\enddemo

The next result deals with the converse implication and is analogous to Proposition \ref{convolutorlemma2}; again we reduce its proof to the weight sequence setting.

\begin{proposition}\label{convolutorlemmaweight2}
Let $r,u,w$ be normalized weights and let $r$ and $u$ be convex.
\begin{itemize}
\item[$(i)$] Assume that
$$\exists\;c_1,c_2,d>0\;\forall\;f\in H^{\infty}_{r_{c_1}}, g\in H^{\infty}_{u_{c_2}}:\;\;\;f\cdot g\in H^{\infty}_{w_{d}}.$$
Then the weights satisfy
$$\exists\;A,B\ge 1\;\forall\;t\ge 0:\;\;\;(w(At))^2\le Br(t)u(t).$$

\item[$(ii)$] If $(i)$ holds for $r\equiv u\equiv w$, then $u$ has to be of moderate growth.
\end{itemize}
Alternatively, in order to conclude in $(i)$ we can replace some or all $\exists$ by $\forall$.
\end{proposition}

\demo{Proof}
$(i)$ Let $c_1,c_2,d>0$ be the parameters. The second estimate in \eqref{omegavequiv} and the assumption imply that for any $f\in H^{\infty}_{v_{M^r},c_1}$, $g\in H^{\infty}_{v_{M^u},c_2}$ we get $f\cdot g\in H^{\infty}_{w_{d}}$.

We apply this to $f\equiv \theta_{M^r,c_1}$ and $g\equiv\theta_{M^u,c_2}$ from Lemma \ref{charholomfctlemma} and restrict to $t\ge 0$ in order to get
$$\exists\;C\ge 1\;\forall\;t\ge 0:\;\;\;\omega_{M^r}(c_1t/2)+\omega_{M^u}(c_2t/2)\le\log(C)-\log(w(dt))=\log(C)+\omega^w(dt),$$
recall \eqref{omegafromv} for the last equality. Set $c:=\min\{c_1,c_2\}$ and thus we have $$\omega_{M^r}(c_1t/2)+\omega_{M^u}(c_2t/2)\ge\omega_{M^r}(ct/2)+\omega_{M^u}(ct/2)=\omega_{M^r\star M^u}(ct/2),$$
see again Remark \ref{convolvesequrem0}. Similarly like in the proof of Proposition \ref{convolutorlemma2} and by taking into account \eqref{vBMTweight1equ1} (for $M^w$) we get for all $j\in\NN$:
\begin{align*}
M^w_j&=\sup_{t>0}\frac{t^j}{\exp(\omega^w(t))}=\sup_{s>0}\frac{(ds)^j}{\exp(\omega^w(ds))}\le Cd^j\sup_{s>0}\frac{s^j}{\exp(\omega_{M^r\star M^u}(cs/2))}
\\&
=C\left(\frac{2d}{c}\right)^j(M^r\star M^u)_j.
\end{align*}
Hence $M^w\hyperlink{preceq}{\preceq}M^r\star M^u$ is verified and by Remarks \ref{convolvesequrem0} and \ref{convolvesequrem} this relation yields
$$\exists\;H\ge 1\;\forall\;t\ge 0:\;\;\;\omega_{M^r\star M^u}(t)=\omega_{M^r}(t)+\omega_{M^u}(t)\le\omega_{M^w}(Ht)+H\Leftrightarrow v_{M^w,H}(t)\le e^Hv_{M^r}(t)v_{M^u}(t).$$
We use this, the first estimate in \eqref{omegavequiv} applied to $r$ and $u$ (for which convexity is required) and the second one applied to $w$ and get
$$\exists\;A\ge 1\;\exists\;H\ge 1\;\forall\;t\ge 0:\;\;\;w(Ht)\le v_{M^w}(Ht)\le e^H v_{M^r}(t)v_{M^u}(t)\le A e^H\sqrt{r(t)u(t)}.$$

$(ii)$ When $r\equiv u\equiv w$, then following the proof in $(i)$ we get $M^{u}\hyperlink{preceq}{\preceq}M^{u}\star M^{u}$ and so $M^{u}$ satisfies \hyperlink{mg}{$(\on{mg})$} by Remark \ref{convolvesequrem}. Thus, by Lemma \ref{om6rem} we get \eqref{om6forv} for $u$.
\qed\enddemo

\begin{proposition}\label{convolutorlemmaweight3}
Let $u$ be a normalized and convex weight function. Assume that $H^{\infty}_{\overline{\mathcal{U}}_{\mathfrak{c}}}$ is closed under point-wise multiplication and that $\mathfrak{m}$ is continuously acting on this space. Then $u$ satisfies \eqref{om6forv}; i.e. $u$ is of moderate growth.
\end{proposition}

\demo{Proof}
Analogously, as in the proof of Proposition \ref{convolutorlemma3} we get $H^{\infty}_{\overline{\mathcal{U}}_{\mathfrak{c}}}\subseteq H^{2,\infty}_{\overline{\mathcal{U}}_{\mathfrak{c}}}$ for the Fr\'{e}chet space
$$H^{2,\infty}_{\overline{\mathcal{U}}_{\mathfrak{c}}}:=\{f\in H(\CC): \|f^2\|_{u_c}<+\infty,\;\forall\; c>0\}.$$
By the closed graph theorem the inclusion $H^{\infty}_{\overline{\mathcal{U}}_{\mathfrak{c}}}\subseteq H^{2,\infty}_{\overline{\mathcal{U}}_{\mathfrak{c}}}$ is continuous and hence
$$\forall\;d>0\;\exists\;c>0\;\exists\;C\ge 1\;\forall\;f\in H^{\infty}_{\overline{\mathcal{U}}_{\mathfrak{c}}}:\;\;\;\|f^2\|_{u_d}\le C\|f\|_{u_c}.$$
We apply this estimate to the monomials $f_k=z^k$ and get (see \eqref{omegafromv}):
\begin{align*}
&\forall\;d>0\;\exists\;c>0\;\exists\;C\ge 1\;\forall\;k\in\NN\;\forall\;z\in\CC:\
\\&
\frac{|z|^{2k}}{\exp(\omega^u(d|z|))}=|z|^{2k}u(d|z|)\le C|z|^ku(c|z|)=C\frac{|z|^k}{\exp(\omega^u(c|z|))}.
\end{align*}
Then, we use \eqref{vBMTweight1equ1} and follow the proof of Proposition \ref{convolutorlemma3} to get that $M^u$ satisfies \eqref{mgdiag} and thus \hyperlink{mg}{$(\on{mg})$}. Consequently, by Lemma \ref{om6rem} we get \eqref{om6forv} for $u$.
\qed\enddemo

When combining Lemma \ref{om6rem} and Propositions \ref{convolutorlemmaweight1}, \ref{convolutorlemmaweight2} and \ref{convolutorlemmaweight3} we arrive at the following characterization.

\begin{theorem}\label{convolutorthmweightfct}
Let $u$ be a normalized and convex weight function. Then the following are equivalent:
\begin{itemize}
\item[$(i)$] $H^{\infty}_{\underline{\mathcal{U}}_{\mathfrak{c}}}$ is closed under point-wise multiplication and $\mathfrak{m}$ is continuously acting on this space.

\item[$(ii)$] $H^{\infty}_{\overline{\mathcal{U}}_{\mathfrak{c}}}$ is closed under point-wise multiplication and $\mathfrak{m}$ is continuously acting on this space.

\item[$(iii)$] $u$ is of moderate growth; i.e. $u$ satisfies \eqref{om6forv}.

\item[$(iv)$] The associated sequence $M^u$ satisfies \hyperlink{mg}{$(\on{mg})$}.
\end{itemize}
\end{theorem}

We close by showing that the proofs of Proposition \ref{convolutorlemma3} and \ref{convolutorlemmaweight3} give an independent argument for the conclusion of Proposition \ref{Pepecounterproof}; even under more general assumptions on the weight.

\begin{corollary}\label{nonpointwiseclosed}
Let $M\in\hyperlink{LCset}{\mathcal{LC}}$ resp. a normalized weight $u$ be given. Then $H^{\infty}_{v_M}$ resp. $H^{\infty}_{u}$ is not closed under point-wise multiplication.
\end{corollary}

\demo{Proof}
For the weight sequence case we present two proofs: First, by Theorem \ref{essentialweightsequthm} we know that $v_M$ is essential and thus Proposition \ref{Pepecounterproof} yields the conclusion.

Second, given $M\in\hyperlink{LCset}{\mathcal{LC}}$ let us assume that $H^{\infty}_{v_M}$ is closed under point-wise multiplication. Then follow the arguments in the proof of Proposition \ref{convolutorlemma3} with $d=c=1$ and thus get
\begin{equation}\label{toostrongformg}
\exists\;C\ge 1\;\forall\;j\in\NN:\;\;\;M_{2j}\le CM_j.
\end{equation}
But this violates $\lim_{j\rightarrow+\infty}(M_j)^{1/j}=+\infty$: Let $k\in\NN$ be such that $2^n\le k<2^{n+1}$ for some $n\in\NN_{>0}$, then $M_k\le M_{2^{n+1}}\le C^{n+1}M_1$ and so $(M_k)^{1/k}\le M_1^{1/k}C^{(n+1)/k}\le M_1^{1/k}C^{(n+1)/2^n}\le M_1C<+\infty$.

Now given a normalized weight $u$ and assuming that $H^{\infty}_{u}$ is closed under point-wise multiplication, then the proof of Proposition \ref{convolutorlemmaweight3} with $d=c=1$ gives \eqref{toostrongformg} for $M^u$. However, since $M^u\in\hyperlink{LCset}{\mathcal{LC}}$ this is impossible by the arguments just given before.
\qed\enddemo

Note:

\begin{itemize}
\item[$(*)$] For the conclusion convexity for $u$ is not required since in Proposition \ref{convolutorlemmaweight3} convexity was only needed in order to apply Lemma \ref{om6rem} in the very last step which is not needed in Corollary \ref{nonpointwiseclosed}.

\item[$(*)$] For the contradictory argument we have used that $\lim_{j\rightarrow+\infty}(M^u_j)^{1/j}=+\infty$. In view of Remark \ref{assoweightsequ} we know that this argument would break down if $u(t)=0$ for all $t$ sufficiently large. But then as sets $H_u^{\infty}=H(\CC)$ and $H(\CC)$ is clearly closed under point-wise multiplication.

\item[$(*)$] The difference between the failure in the single weight case and the characterization in the weight system case is that the latter situation allows for having ''enough place'' when changing the parameter; see the proof of Proposition \ref{poinwisepowerclear} and note that $H>1$ in \eqref{om6forv}.

\item[$(*)$] The characterization \cite[Prop. 2.1]{BonetDomanskiLindstroem99} suggests that Corollary \ref{nonpointwiseclosed} fails for classes $H^{\infty}_{v_M}(\mathbb{D})$ resp. $H^{\infty}_{u}(\mathbb{D})$. Note that by \cite[Sect. 7]{solidassociatedweight}, in the weight sequence setting, for $G=\mathbb{D}$ the weights are given by $t\mapsto\exp(-\omega_M(\frac{1}{1-t}))=:v_{M,\mathbb{D}}(t)$ and hence the corresponding norms are expressed in terms of $v_{M,\mathbb{D}}$. Consequently, when following the proof of Proposition \ref{convolutorlemma3} we arrive in \eqref{convolutorlemma3equ} (again with $c=d=1$) at expressions which are not well-behaved w.r.t. the crucial formula \eqref{Prop32Komatsu}.
\end{itemize}

\bibliographystyle{plain}
\bibliography{Bibliography}

\end{document}